\documentclass[12pt,leqno]{article}

\usepackage{times}
\usepackage{a4wide}
\usepackage{amssymb}
\usepackage{amsmath}
\usepackage{theorem}
\usepackage{enumerate}
\usepackage{amsmath}
\usepackage{amscd}
\usepackage{amssymb}

\theoremstyle{change}
\theorembodyfont{\sl}
\newtheorem{theorem}[subsection]{Theorem.}
\newtheorem{proposition}[subsection]{Proposition.}
\newtheorem{lemma}[subsection]{Lemma.}
\newtheorem{corollary}[subsection]{Corollary.}

\newtheorem{sublemma}[subsubsection]{Lemma.}

\theorembodyfont{\rmfamily}
\newtheorem{remark}[subsection]{Remark.}
\newtheorem{subremark}[subsubsection]{Remark.}

\makeatletter
\newenvironment{eqn}{\refstepcounter{subsection}
$$}{\leqno{\rm(\thesubsection)}$$\global\@ignoretrue}
\makeatother

\makeatletter 
\newenvironment{subeqn}{\refstepcounter{subsubsection}
$$}{\leqno{\rm(\thesubsubsection)}$$\global\@ignoretrue}
\makeatother 



\newenvironment{figuur}%
{%
\refstepcounter{subsection}\begin{figure}}%
{\end{figure}}

\newenvironment{prf}[1]{\trivlist
\item[\hskip \labelsep{\it
#1.\hspace*{.3em}}]}{~\hspace{\fill}~$\square$\endtrivlist}

\newenvironment{proof}{\begin{prf}{\bf Proof}}{\end{prf}}

\newcommand{\ol}{\overline}
\newcommand{\ZZ}{{\mathbb Z}}
\newcommand{\NN}{{\mathbb N}}
\newcommand{\QQ}{{\mathbb Q}}

\newcommand{\CC}{{\mathbb C}}
\newcommand{\FF}{{\mathbb F}}

\newcommand{\PP}{{\mathbb P}}
\newcommand{\GG}{{\mathbb G}}

\renewcommand{\AA}{{\mathbb A}}

\newcommand{\Pic}{{\rm Pic}}
\newcommand{\Spec}{{\rm Spec}}
\newcommand{\rH}{{\rm H}}

\newcommand{\calO}{{\cal O}}
\newcommand{\calL}{{\cal L}}
\newcommand{\calM}{{\cal M}}
\newcommand{\calC}{{\cal C}}
\newcommand{\calZ}{{\cal Z}}
\newcommand{\calMbar}{\ol{\calM}}
\newcommand{\calCbar}{\ol{\calC}}
\newcommand{\swedge}{{\scriptstyle{\wedge}}}
\newcommand{\nor}{{\rm nor}}
\newcommand{\divisor}{{\rm div}}
\newcommand{\rmm}{{\rm m}}
\newcommand{\Div}{{\rm Div}}
\newcommand{\Sing}{{\rm Sing}}
\newcommand{\Norm}{{\rm Norm}}

\newcommand{\GL}{{\rm GL}}
\newcommand{\om}{\underline{\omega}}
\newcommand{\cusps}{{\rm cusps}}
\newcommand{\rmD}{{\rm D}}
\newcommand{\lto}{\longrightarrow}
\newcommand{\Hom}{{\rm Hom}}
\newcommand{\coker}{{\rm coker}}
\newcommand{\im}{{\rm im}}

\begin{document}
\title{On N\'eron models, divisors and modular curves.} 
\author{Bas Edixhoven\footnote{partially supported by the 
Institut Universitaire de France, and by the European TMR Network 
Contract ERB FMRX 960006 ``arithmetic algebraic geometry''.}}
\maketitle

\section{Introduction.}\label{section1}
For $p$ a prime number, let $X_0(p)_\QQ$ be the modular curve, over
$\QQ$, parametrizing isogenies of degree $p$ between elliptic curves,
and let $J_0(p)_\QQ$ be its jacobian variety. Let $f\colon
X_0(p)_\QQ\to J_0(p)_\QQ$ be the morphism of varieties over $\QQ$ that
sends a point $P$ to the class of the divisor $P-\infty$, where
$\infty$ is the $\QQ$-valued point of $X_0(p)_\QQ$ that corresponds to
$\PP^1_\QQ$ with $0$ and $\infty$ identified, equipped with the
subgroup scheme~$\mu_{p,\QQ}$.  We suppose that $X_0(p)_\QQ$ has genus
at least one, i.e., that $p=11$ or $p>13$. Then $f$ is a closed
immersion.  Let $0$ denote the other cusp of~$X_0(p)_\QQ$. The
$\QQ$-valued point $f(0)$ of $J_0(p)_\QQ$ is well known to be of order
$n$, the numerator of $(p-1)/12$, expressed in lowest terms. Let now
$X_0(p)$ denote the model over $\ZZ$ of $X_0(p)_\QQ$ as described by
\cite{DeligneRapoport1} and \cite{KatzMazur1}: it is the compactified
coarse moduli space for generalized elliptic curves with finite
locally free subgroup schemes of rank $p$ that meet all irreducible
components of all geometric fibres.  Let $X_0(p)^\sim$ denote the
minimal regular model of $X_0(p)$, and let $J_0(p)$ be the N\'eron
model over $\ZZ$ of $J_0(p)_\QQ$.  See \cite{DeligneRapoport1} for a
description of the semi-stable curve~$X_0(p)^\sim$. By the defining
properties of $J_0(p)$, the morphism $f$ extends uniquely to a
morphism $f\colon X_0(p)^{\sim,0}\to J_0(p)$, where $X_0(p)^{\sim,0}$
is the open part of $X_0(p)^\sim$ where the morphism to $\Spec(\ZZ)$
is smooth; $X_0(p)^{\sim,0}$ is the complement of the set of double
points in the fibre $X_0(p)^\sim_p$ over~$\FF_p$. Robert Coleman asked
about how the image of $X_0(p)^{\sim,0}_p$ under $f$ intersects $C$,
the closure in $J_0(p)$ of the group generated by~$f(0)$. We will
prove, see Theorem~\ref{theorem8.2}, for all $p$ for which
$X_0(p)_\QQ$ has genus at least two, that the intersection consists
just of the two obvious elements $f(0)_p$ and~$f(\infty)_p$. (Of
course, when $X_0(p)_\QQ$ has genus one, $f$ is an isomorphism, hence
the intersection is all of~$C_p$.) This result is used by Coleman,
Kaskel and Ribet in \cite{ColKasRib} to study the inverse image under
$f$ of the torsion subgroup of~$J_0(p)(\CC)$.

To prove our result, we use Raynaud's description (see \cite{Raynaud1}
or \cite[9.5]{BLR}) of the N\'eron model of a jacobian variety of a
curve in terms of a regular model of the curve. (A good reference for
jacobians of modular curves is~\cite{Raynaud2}.)  It turns out that in
order to describe the special fibre of such a N\'eron model, and not
just the connected component of the identity element and the group of
connected components, one needs to know a bit more than just the
special fibre of the regular model of the curve. This extra
information is related, in the semi-stable case, to the log structure
induced on the special fibre of the curve (see
Section~\ref{section6}), and also to the map from the space of global
deformations to that of local deformations of a stable curve
(see~\cite[p.~81]{DeligneMumford}).  For the curve $X_0(p)^\sim$ we
describe this extra information in terms of a well known modular form
of weight $p+1$ on the supersingular elliptic curves. We do not use
rigid uniformization, and our methods can be used as well for curves
whose special fibre is not totally degenerate.  A rigid analytic
description of $J_0(p)_p$ is given by de Shalit in~\cite{deShalit1}.

In Sections~\ref{section2} and~\ref{section3} we consider a proper
nodal curve $X$ with smooth generic fibre $X_K$ over a discrete
valuation ring, and study the open part $J$ of the N\'eron model of
the jacobian variety given by the relative Picard functor of~$X$. We
give a description of the special fibre $J_k$ of $J$ in terms of
divisors on $X_k$ and some invariants $c(x)$ at the double points
of~$X_k$. Some properties of these $c(x)$ are discussed in
Sections~\ref{section4}, \ref{section5} and \ref{section6}. The $c(x)$
for the modular curve $X_0(p)$ are determined in
Section~\ref{section7}.  Coleman's question mentioned above is
answered in Section~\ref{section8}.  Section~\ref{section9} answers
the question as to what extent the morphism from the smooth locus
$X^0$ of $X$ to the N\'eron model of $J_K$ is a closed immersion. For
example, if $X$ is regular, we show that this morphism is a closed
immersion if and only if all double points of $X_{\bar{k}}$ are
non-disconnecting (see Corollary~\ref{cor9.8}).  Some of the results
in this last section can also be found in Coleman's
preprint~\cite{Coleman}.

\section{Raynaud's description in the nodal case.}\label{section2}
Let $D$ be a complete discrete valuation ring, let $\pi$ be a
uniformizer, $K$ the fraction field of $D$ and $k$ the residue
field. We suppose that $k$ is separably closed. Let $X_K$ be a smooth,
proper and geometrically irreducible curve over $K$, and let $X$ be a
proper flat model of it over~$D$. We suppose that $X$ is nodal, i.e.,
the only singularities of the geometric special fibre are ordinary
double points.  The irreducible components of $X_k$ are absolutely
irreducible.  The singular points of $X_k$ are $k$-rational points. At
such a point $x$ the complete local ring $\calO_{X,x}^\swedge$ is
isomorphic to one of the form $D[[u,v]]/(uv-\pi^{e(x)})$, for a unique
integer $e(x)\geq1$.  Let $R\colon X^\sim\to X$ be the minimal
resolution of singularities of~$X$; for $x$ in $X_k$ a singular point,
$R^{-1}x$ is a chain of $e(x)-1$ projective lines.

Let $N\colon X_k^\nor\to X_k$ be the normalization morphism. We define
$S_0$ to be the set of irreducible components of $X_k^\nor$. Note that
$N$ induces a bijection between $S_0$ and the set of irreducible
components of~$X_k$. Let $S_1$ be the set of singular points of $X_k$,
and let $S_1^\nor$ be~$N^{-1}S_1$. We will view $S_0$ and $S_1$ as the
sets of vertices and edges of a graph; the end points of $x$ in $S_1$
are the elements of $S_0$ containing at least one of the two elements
of~$N^{-1}\{x\}$.

According to \cite[Thm.~8.2.1]{Raynaud1} or \cite[Thm.~8.2/2]{BLR},
the relative Picard functor $\Pic_{X/D}$ is representable, and is
smooth over $D$ (but not of finite type, and not separated if $X_k$
has more than one irreducible component).  As usual, let
$\Pic_{X/D}^0$ be the connected component of the identity element in
$\Pic_{X/D}$; its $D$-valued points correspond to line bundles on $X$
whose restrictions to all irreducible components of $X_k$ are of
degree zero. As in \cite{Mazur1}, let $\Pic_{X/D}^{[0]}$ be the
closure in $\Pic_{X/D}$ of $\Pic_{X_K/K}^0$; $\Pic_{X/D}^{[0]}(D)$ is
the set of isomorphism classes of line bundles on $X$ of which the sum
of the degrees of the restrictions to the irreducible components is
zero. Let $E$ be the closure in $\Pic_{X/D}$ of the identity element
of~$\Pic_{X_K/K}$. Then $E(D)$ is the set of isomorphism classes of
line bundles on $X$ whose restriction to $X_K$ is trivial.  Such line
bundles are given by Cartier divisors with support in~$X_k$.  The
group of Weil divisors with support in $X_k$ is the group~$\ZZ^{S_0}$.
Let $F$ be the subgroup of $\ZZ^{S_0}$ consisting of those Weil
divisors that define Cartier divisors. Of course, if $X$ is regular,
then $F$ equals~$\ZZ^{S_0}$. Let us consider a Weil divisor $a:=\sum
a(C)C$ in~$\ZZ^{S_0}$. Then $a$ defines a Cartier divisor if and only
if $a$ is principal at all $x$ in~$S_1$. Let $x$ be in $S_1$, and let
$C$ and $C'$ in $S_0$ be the endpoints of~$x$. Then $a$ is principal
at $x$ if and only if $a(C)$ and $a(C')$ are equal modulo~$e(x)$ (see
\cite[p.~15]{Raynaud2}).  It follows that we have an exact sequence:
\begin{eqn}\label{eqn2.1}
0 \to F\to\ZZ^{S_0} \to 
\prod_{x\in S_1} (\ZZ/e(x)\ZZ)^{N^{-1}\{x\}}/(\ZZ/e(x)\ZZ), 
\end{eqn}
where each $\ZZ/e(x)\ZZ$ is diagonally embedded in~$\prod_{x\in S_1}
(\ZZ/e(x)\ZZ)^{N^{-1}\{x\}}$. By definition, we have a map from $F$ to
$\Pic_{X/D}^{[0]}(D)$, whose image is~$E(D)$. This map sends $a$ in
$F$ to the invertible $\calO_X$-module $\calO_X(a)$, whose sections on
open $U\subset X$ are rational functions $f$ on $U$ with
$\divisor(f)+a\geq0$.

By \cite{Raynaud1}, or the proof of \cite[Prop.~9.5/3]{BLR}, $E$ is
the espace \'etal\'e of the sheaf $i_*E(D)$, with
$i\colon\Spec(k)\to\Spec(D)$ the closed immersion given by the
canonical projection $D\to k$. In particular, $E(k)=E(D)$, and
$E(K)=0$.  It follows that $J:=\Pic_{X/D}^{[0]}/E$ is representable
and separated, and that $J$ represents the open part of the N\'eron
model over $D$ of $\Pic^0_{X_K/K}$ that corresponds to the line
bundles on $X_K$ that admit an extension to a line bundle on~$X$.
Note that if $X$ is regular, then $J$ is the N\'eron model itself.

Let $\Phi:=J_k/J_k^0$ be the group of connected components of~$J_k$. 
By \cite[Thm.~9.7/1]{BLR}, $J_k^0=\Pic^0_{X_k/k}$. 
We have the following complex of $\ZZ$-modules: 
\begin{eqn}\label{eqn2.2}
0 \to \ZZ \to F\to \ZZ^{S_0} \to \ZZ \to 0. 
\end{eqn}
The element $1$ in $\ZZ$ is mapped to $X_k$, the principal divisor
defined by $\pi$. The morphism from $F$ to $\ZZ^{S_0}$ is the
composite of $F\to\Pic(X)$ and the morphism ``multidegree'' from
$\Pic(X)$ to $\ZZ^{S_0}$ that sends a line bundle to its degrees on
the irreducible components of~$X_k$ (note that this map is not the
same as the one in~(\ref{eqn2.1})). Explicitly, for $a$ in $F$ and $C$
in $S_0$, one has:
\begin{eqn}\label{eqn2.3}
\deg_C(\calO_X(a)) = \sum_{x\in C\cap S_1^\nor}\frac{a(C'_x)-a(C)}{e(x)}, 
\end{eqn}
where $C'_x$ denotes the other branch passing through $x$
(see~(\ref{eqn3.5})).  The morphism $\ZZ^{S_0}\to\ZZ$ is simply the
sum. The fact that the $\QQ$-valued intersection pairing on
$\ZZ^{S_0}$ is negative definite on $\ZZ^{S_0}/\ZZ X_k$ shows that
(\ref{eqn2.2}) is exact at~$F$.  The definition of $J$, and the fact
that $E(D)$ is the image of $F$, show that the homology of
(\ref{eqn2.2}) at $\ZZ^{S_0}$ is~$\Phi$.  See \cite[Thm.~9.6/1]{BLR}
for this result in the case where $X$ is regular.

\section{A description of $J_k$ in terms of divisors.} \label{section3}
We keep the notation and hypotheses of the preceding section.  The
results in that section give us a description of $J$ in terms
of~$\Pic_{X/D}$. Further on, we will have to do computations in the
special fibre $J_k$ of $J$, without knowing too much about the
$D$-scheme $X$, but in terms of the $k$-scheme $X_k$ with some
additional data.

In Section~\ref{section2} we have seen that $J_k$ has the following 
presentation: 
\begin{eqn}\label{eqn3.1}
F \to \Pic^{[0]}_{X_k/k} \to J_k \to 0, 
\end{eqn}
with $\Pic^{[0]}_{X_k/k}$ the open part of $\Pic_{X_k/k}$
corresponding to line bundles on $X_k$ with total degree zero, and
with $F$ the subgroup of $\ZZ^{S_0}$ consisting of the Weil divisors
that are locally principal (we view $F$ as a constant group scheme
over~$k$).  This subgroup $F$ is described in terms of $S_0$ and the
$e(x)$ in~(\ref{eqn2.1}). The only ingredient in (\ref{eqn3.1}) that
does not only involve $X_k$ and the $e(x)$ is the map from $F$
to~$\Pic(X_k)$. The aim of this section is to describe $\Pic_{X_k/k}$
and the map from $F$ to it, in terms of $X_k$, the $e(x)$, and some
extra data related to the reduction of $X$ modulo some power of $\pi$
(recall that $\pi$ is a uniformizer of~$D$).

Let us first describe $\Pic_{X_k/k}$. The proof of Example~8 of
Chapter~9 of \cite{BLR} shows that we have an exact sequence:
\begin{eqn}\label{eqn3.2}
0\to \GG_{\rmm,k} \to \bigoplus_{C\in S_0}\GG_{\rmm,k} \to 
\bigoplus_{x\in S_1}\left(\GG_{\rmm,k}^{N^{-1}\{x\}}/\GG_{\rmm,k}\right) 
\to \Pic_{X_k/k} \to \Pic_{X_k^\nor/k} \to 0. 
\end{eqn}
Let us assume that $S_1$ is not empty. Then (\ref{eqn3.2}) gives the
following description of $\Pic_{X_k/k}(S)$ for all $k$-schemes~$S$.
\begin{proposition}\label{prop3.3}
Under the assumptions above, $\Pic_{X_k/k}(S)$ is the set of
isomorphism classes of line bundles on $X_S$ that are obtained as
follows.  Choose a line bundle $\calL$ on $X^\nor_S$, i.e., for each
$C$ in $S_0$, a line bundle $\calL_C$ on $C_S$, with a given
isomorphism $\phi_x\colon \calO_S\to x^*\calL$ for each $x$
in~$S_1^\nor$. Then glue, at each $x$ in $S_1$, the line bundles
$\calL_{C_1}$ and $\calL_{C_2}$ via $\phi_{x_1}$ and $\phi_{x_2}$,
where $\{x_1,x_2\}$ is $N^{-1}\{x\}$ and where $C_1$, $C_2$ are the
elements of $S_0$ with $x_1\in C_1$ and $x_2\in C_2$. Let
$(\calL,\phi)$ denote the line bundle on $X_S$ obtained in this
way. Two such line bundles $(\calL,\phi)$ and $(\calL',\phi')$ are
isomorphic if and only if there exists an isomorphism
$\alpha\colon\calL\to\calL'$ that is compatible with the glueing data
in the sense that for all $x$ in $S_1$ we have
$$
\alpha_{x_2}\circ\phi_{x_2}\circ\phi_{x_1}^{-1} = 
\phi'_{x_2}\circ{\phi'_{x_1}}^{-1}\circ\alpha_{x_1}.
$$
\end{proposition}

\begin{remark}\label{remark3.4}
Another way to state the previous proposition is to say that 
$N\colon X_k^\nor\to X_k$ is a morphism of universal effective descent 
for line bundles (see \cite[D\'ef.~1.7]{Grothendieck1}). 
\end{remark}
Our next aim is to describe the line bundles of the form
$\calO_X(a)|_{X_k}$, with $a$ in $F$, in the form of
Proposition~\ref{prop3.3}. Note that the divisor of $\pi$ is
$X_k=\sum_CC$, so that multiplication by $\pi$ induces an isomorphism
from $\calO_X$ to~$\calO_X(-X_k)$. Let $a$ be in $F$.  For $C$ in
$S_0$ we have:
\begin{eqn}\label{eqn3.5}
\calO_X(a)|_C \;\;\tilde{\to}\;\; 
\left.\calO_X\left(\sum_{C'\neq C}(a(C')-a(C))C'\right)\right|_C = 
\calO_C\left(\sum_{x\in C\cap S_1^\nor}\frac{a(C'_x)-a(C)}{e(x)}x\right), 
\end{eqn}
where the first isomorphism is given by multiplication by
$\pi^{a(C)}$, and where $C'_x$ is the other branch passing
through~$x$.  Now we have to see how these line bundles are glued at
the singular points. So let $x$ be in $S_1$,
$N^{-1}\{x\}=\{x_1,x_2\}$, and let $C_1$ and $C_2$ be the elements of
$S_0$ containing $x_1$ and $x_2$, respectively. Let $a_1:=a(C_1)$ and
$a_2:=a(C_2)$. Then we have:
\begin{eqn}\label{eqn3.6}
x_1^*\calO_{C_1}\left(\frac{a_2-a_1}{e(x)}x_1\right) = 
\Omega_{C_1}(x_1)^{\otimes\frac{a_1-a_2}{e(x)}}, 
\end{eqn}
where $\Omega_{C_1}$ is the sheaf of relative 1-forms of $C_1$
over~$k$.  Of course, we have a similar formula at~$x_2$. It follows
that multiplication by $\pi^{a_1-a_2}$ defines an isomorphism:
\begin{eqn}\label{eqn3.7}
\Omega_{C_2}(x_2)^{\otimes\frac{a_2-a_1}{e(x)}} \to 
\Omega_{C_1}(x_1)^{\otimes\frac{a_1-a_2}{e(x)}}, 
\end{eqn}
or, equivalently, a non-zero element of
$(\Omega_{C_1}(x_1)\otimes\Omega_{C_2}(x_2))^{\otimes(a_2-a_1)/e(x)}$.
A simple computation shows that if we choose an isomorphism
$D[[u,v]]/(uv-\pi^{e(x)})\to\calO_{X,x}^\swedge$, such that $u$ is
zero on $C_1$ and $v$ is zero on $C_2$, then this non-zero element is
simply the $(a_2-a_1)/e(x)$th power of
\begin{eqn}\label{eqn3.8}
\mbox{$c(x):=dv\otimes du$ in $\Omega_{C_1}(x_1)\otimes\Omega_{C_2}(x_2)$} 
\end{eqn}
(we will see in the next section that $c(x)$ does not depend on the
choice of the isomorphism).  This finishes our description of the map
from $F$ to $\Pic(X_k)$, and hence of the description of~$J_k$. For
the sake of notation, let us say that for $a$ in $F$, we call
$\calL_a$ the line bundle on $X_k$ that we just constructed (so
$\calL_a$ is isomorphic to $\calO_X(a)|_{X_k}$, but its description in
terms of Proposition~\ref{prop3.3} depends on $\pi$).  Some properties
of the $c(x)$ will be discussed in the next section. To finish this
section, we describe $J_k$ in terms of degree zero divisors on~$X_k$.

Let $\Div(X_k)$ be the $\ZZ$-module of Weil divisors on $X_k$ with
support outside the singular locus~$S_1$. Since such divisors are
locally principal, we have a morphism $\Div(X_k)\to\Pic(X_k)$, sending
$E$ in $\Div(X_k)$ to isomorphism class of the line
bundle~$\calO_{X_k}(E)$.  This morphism is easily seen to be
surjective. Let $\Div^{[0]}(X_k)$ be the subgroup of $\Div(X_k)$
consisting of the divisors of degree zero. Recall that $J_k$ is the
quotient of $\Pic^{[0]}_{X_k/k}$ by the image of~$F$.  We get a
morphism $\Div^{[0]}(X_k)\to J_k(k)$, sending $E$ to $P_E$, say.  It
remains to describe the kernel of $E\mapsto P_E$.

So let $E$ be in $\Div^{[0]}(X_k)$. Of course, $P_E$ is zero if and
only if there exists an $a$ in $F$ such that $\calL_a$ is isomorphic
to~$\calO_{X_k}(E)$. A necessary condition for $P_E$ to be zero is
that the image of $P_E$ in $\Phi$ be zero. This last condition is
equivalent to the multidegree of $E$ in $\ZZ^{S_0}$ being in the image
of $F$ (see (\ref{eqn2.2})); let us suppose now that that is the case.
Let $a=\sum_Ca(C)C$ be in $F$ such that $\calL_a$ and $E$ have the
same multidegree; $a$ is unique up to adding multiples of~$X_k$. Then
$P_E$ is zero if and only if $\calL_a\otimes\calO_{X_k}(-E)$ is
trivial, or, equivalently, if $\calL_a\otimes\calO_{X_k}(-E)$ has a
non-zero global section (note that it has degree zero on
each~$C$). Using the description in Proposition~\ref{prop3.3} of line
bundles on $X_k$, the last condition amounts to saying that the
$\calL_{a,C}\otimes\calO_C(-E|_C)$ have non-zero global sections $f_C$
that are compatible at all $x$ in~$S_1$.  By construction,
see~(\ref{eqn3.5}), each such $f_C$ has to be a non-zero rational
function on $C$, whose divisor is:
\begin{eqn}\label{eqn3.9}
E|_C - \sum_{x\in C\cap S_1^\nor}\frac{a(C'_x)-a(C)}{e(x)}x. 
\end{eqn}
Suppose now that there exist $f_C$ with these divisors, or,
equivalently, suppose that the $\calL_{a,C}\otimes\calO_C(-E|_C)$ are
trivial. Let $f_C$ be such functions. In order to make the
compatibility conditions for the $f_C$ explicit, we introduce some
terminology.  For $f$ a non-zero rational function on $C$ in $S_0$,
and $x$ in $C(k)$, we define the leading term $f(x)$ of $f$ at $x$ to
be the non-zero element of $\Omega_C(x)^{\otimes n}=m_x^n/m_x^{n+1}$
given by $f$, where $n$ is the order of $f$ at $x$, and $m_x$ the
maximal ideal in~$\calO_{C,x}$.  For $x$ in $S_1$, $C_1$ and $C_2$ the
elements of $S_0$ containing the preimages $x_1$ and $x_2$ of $x$
under $N\colon X_k^\nor\to X_k$, the two rational functions $f_{C_1}$
and $f_{C_2}$ are compatible at $x$ if and only if:
\begin{eqn}\label{eqn3.10}
f_{C_1}(x_1)\,f_{C_2}(x_2)^{-1}c(x)^{(a(C_2)-a(C_1))/e(x)} = 1, 
\end{eqn}
as can be seen from the definitions earlier in this section. The $f_C$
are unique up to multiplication by elements from~$k^*$. The problem of
whether or not one can scale the $f_C$ such that they satisfy the
compatibility conditions at all $x$ in $S_1$ is a problem concerning
the structure of the graph $(S_0,S_1)$ of~$X_k$.  The following
proposition gives the final conclusion of the preceding discussion.
\begin{proposition}\label{prop3.11}
Let $E$ be in $\Div^{[0]}(X_k)$. Then $P_E$ in $J_k(k)$ is zero if and
only if the following conditions are satisfied.  First of all, the
multidegree of $E$ in $\ZZ^{S_0}$ has to be in the image of $F$;
suppose this is so and let $a$ be in $F$ having the same image in
$\ZZ^{S_0}$ as~$E$. Secondly, for each $C$ in $S_0$, the divisor
in~(\ref{eqn3.9}) has to be a principal divisor; suppose that this is
so, and let, for each $C$, $f_C$ be a rational function on $C$ with
that divisor. Thirdly, for every cycle
$(C_0,x_0,C_1,x_1,\ldots,x_{m-1},C_m=C_0)$ in the graph $(S_0,S_1)$ we
must have:
\begin{eqn}\label{eqn3.12}
\prod_{i=0}^{m-1}
f_{C_i}(x_i')\,f_{C_{i+1}}(x_i'')^{-1}c(x_i)^{(a(C_{i+1})-a(C_i))/e(x_i)} 
= 1, 
\end{eqn}
with $x_i'$ and $x_i''$ in $N^{-1}\{x\}$ defined by $x_i'\in C_i$ and 
$x_i''\in C_{i+1}$. 
\end{proposition}
\begin{remark}\label{rem3.13}
Of course, in the third condition in Proposition~\ref{prop3.11}, it 
is sufficient to consider a set of cycles that generates the first 
homology group of the graph. 
\end{remark}

\section{Some properties of the $c(x)$.}\label{section4}
Let us briefly recall the definition, given in (\ref{eqn3.8}), of
the~$c(x)$. In the situation of (\ref{eqn3.8}), $D$ is a complete
discrete valuation ring with separably closed residue field $k$ and
fraction field $K$, $\pi$ in $D$ is a uniformizer and $X/D$ is a
proper nodal curve over $D$ with $X_K$ smooth and geometrically
irreducible.  Suppose that $x$ is a singular point of~$X_k$. Then the
$D$-algebra $\calO_{X,x}^\swedge$ is of the type
$D[[u,v]]/(uv-\pi^{e(x)})$ for some unique integer $e(x)>0$. Let $x_1$
and $x_2$ be the two points of $X_k^\nor$ lying over $x$, and let
$C_1$ and $C_2$ be the irreducible components of $X_k^\nor$ that they
lie on.  Let $u$ and $v$ be elements in $\calO_{X,x}^\swedge$ that
vanish on $C_1$ and $C_2$, respectively, that induce parameters, still
denoted $u$ and $v$, of $\calO_{C_2,x_2}^\swedge$ and
$\calO_{C_1,x_1}^\swedge$, respectively, and such that
$uv=\pi^{e(x)}$.  Then we claim that the non-zero element $dv\otimes
du$ in $\Omega_{C_1}(x_1)\otimes\Omega_{C_2}(x_2)$ does not depend on
the choice of $u$ and~$v$. To see this, suppose that $u'$, $v'$ also
satisfy the conditions. Then $u'=\lambda u$ and $v'=\lambda^{-1}v$
with $\lambda$ in $\calO_{X,x}^\swedge$ a unit. Let $\ol{\lambda}$ be
the image in $k^*$ of~$\lambda$. Then $du'=\ol{\lambda}du$, and
$dv'=\ol{\lambda}^{-1}dv$, proving our claim.

Hence the elements $c(x)$ of (\ref{eqn3.8}) are independent of the
choice of the isomorphism. Let us now see how they depend on the
choice of~$\pi$. A simple computation shows that if $\pi'=\lambda\pi$,
then $c'(x)=\ol{\lambda}^{e(x)}c(x)$, where $c'(x)$ is defined
using~$\pi'$.  It follows that the construction of $c(x)$ defines an
isomorphism of $k$-vector spaces between
$\Omega_{C_1}(x_1)\otimes\Omega_{C_2}(x_2)$ and $(m_D/m_D^2)^{\otimes
e(x)}$, where $m_D$ is the maximal ideal of~$D$.

The construction of $c(x)$ makes it clear that $c(x)$ depends only on
the reduction of the $D$-scheme $X$ mod~$\pi^{e(x)+1}$. It follows
that $J_k$, as defined in Section~\ref{section2}, depends only on $X$
mod $\pi^{e+1}$, where $e$ is the maximum of all the~$e(x)$. In
particular, if $X$ is regular, the $c(x)$ and $J_k$ depend only on $X$
mod~$\pi^2$.

The construction of the $c(x)$ also makes sense in more general
situations (which we won't use in this article).  For example,
consider a semi-stable curve $f\colon C\to S$ in the sense of
\cite[2.21]{deJong1}, i.e., $S$ is an arbitrary scheme and $f$ is of
finite presentation, flat, proper and all its geometric fibres are
connected curves whose singularities are ordinary double points. Let
$Z:=\Sing(f)$ be the closed subscheme of $C$ defined by the first
Fitting ideal of $\Omega_{C/S}$, and let $g\colon Z\to S$ be the
inclusion, followed by~$f$. As explained in
\cite[2.21--2.23]{deJong1}, $g$ is of finite presentation, finite and
unramified, and the kernel $I$ of $g^{-1}\calO_S\to\calO_Z$ is locally
principal. We assume that $I/I^2$ is a faithful
$g^{-1}\calO_S/I$-module (this is the case, for example, if $S$ is
integral and $f$ generically smooth).  Then
$\calO_Z\otimes_{g^{-1}\calO_S}I$ is an invertible $\calO_Z$-module.
Let $P$ be the tautological section of $f_Z\colon C_Z\to Z$, and let
$b\colon C_Z^\sim\to C_Z$ be the blow up of $C_Z$ in the ideal sheaf
defined by~$P$. Let $Z'$ be the fibered product of $P$ and $b$, and
denote $P'$ the morphism from $Z'$ to~$C_Z^\sim$. Then $Z'\to Z$ is
the etale $\ZZ/2\ZZ$-torsor given by the two branches of $C_Z$ at $P$,
and $C_Z^\sim\to Z$ is smooth at~$P'$. The global version of the
construction of the $c(x)$ is then an isomorphism of invertible
$\calO_Z$-modules:
\begin{eqn}\label{eqn4.1}
\Norm_{Z'/Z}\left({P'}^*\Omega_{C_Z^\sim/Z}\right) \longrightarrow 
\calO_Z\otimes_{g^{-1}\calO_S}I. 
\end{eqn}
To prove all this, one reduces to the noetherian case and uses the 
description of \cite[2.23]{deJong1}. 

We give a more geometric description of the global version of the
$c(x)$ for the universal stable curve of some genus $g\geq2$. So let
$\calMbar$ be the algebraic stack over $\ZZ$ classifying such curves,
and let $f\colon \calCbar\to\calMbar$ be the universal curve (see
\cite{DeligneMumford} and~\cite{Mumford}).  Put
$\calZ:=\Sing(f)$. Then $\calZ$ is the disjoint union of irreducible
smooth closed substacks $\Delta_i^*$, $0\leq i\leq \lfloor
g/2\rfloor$, of $\calCbar$, of codimension two, and the morphism
$\calZ\to\calMbar$ is locally a closed immersion, defined by a locally
principal ideal $I$ in $f^{-1}\calO_{\calMbar}$. Then then conormal
bundle $I/I^2$ of the locally closed immersion $f\colon\calZ\to\calM$
is a line bundle on~$\calZ$.  There is another line bundle, say
$\calL$, on $\calZ$, with the property that for $x$ in $\calZ(k)$
corresponding to a pair $(C,P)$ with $C$ a stable curve of genus $g$
over the algebraically closed field $k$ and $P$ a singular point of
$C$, $\calL(x)=\Omega_{C_1}(x)\otimes\Omega_{C_2}(x)$, with $C_1$ and
$C_2$ the two branches at~$x$. Then the global version of the $c(x)$
gives an isomorphism between $\calL$ and~$I/I^2$.  The existence of
this isomorphism was already noted, over $\CC$, in
\cite[p.~143]{Looijenga} and in \cite[p.~477]{Pikaart}.

\section{The $c(x)$ and certain morphisms.}\label{section5}
Let $X$, $D$, etc., be as in the beginning of Section~\ref{section2}.
In particular, $R\colon X^\sim\to X$ is the minimal resolution of the
singularities of~$X$. We want to compare the $c(x)$ on $X$
and~$X^\sim$.  So let $x$ be singular in $X(k)$, and suppose we have
an isomorphism between $\calO_{X,x}^\swedge$ and
$D[[u,v]]/(uv-\pi^e)$. Let $C_1$ and $C_2$ be the irreducible
components of $X^\nor$ on which $u$ and $v$ vanish, respectively. Then
$R^{-1}\{x\}$ is the union of $e-1$ smooth, proper and geometrically
irreducible curves of genus zero, say $L_1,\ldots,L_{e-1}$, arranged
in a chain: if we let $L_0:=C_1$ and $L_e:=C_2$, then $L_i\cap L_j$ is
empty if $|i-j|>1$, and $L_{i-1}\cap L_i$ is a unique point $P_i$ for
$1\leq i\leq e$.  Let $u_i:=\pi^{i+1-e}u$ and $v_i:=\pi^{-i}v$.  Then
$\calO_{X^\sim,P_i}^\swedge$ is $D[[u_i,v_i]]/(u_iv_i-\pi)$, $v_i$ is
a parameter of $L_{i-1}$ at $P_i$ and $u_i$ is a parameter of $L_i$
at~$P_i$.  So we have $c(P_i)=dv_i\otimes du_i$. In order to relate
the $c(P_i)$ to $c(x)$ we use the following lemma, the proof of which
is obvious.
\begin{lemma}\label{lem5.1}
On a two-pointed smooth projective curve of genus zero $(L,P_1,P_2)$
over a field $k$, there is a unique non-zero element $x$ in
$\Omega_L(P_1)\otimes\Omega_L(P_2)$ with the following property: let
$f$ be a rational function on $L$ with divisor $P_1-P_2$, then
$x=(df)P_1\otimes (df^{-1})P_2$, with $(df)P_1$ and $df^{-1}P_2$ the
values of $df$ and $df^{-1}$ at $P_1$ and $P_2$, respectively. This
element does not depend on the numbering of the two points.
\end{lemma}
This lemma gives us an isomorphism 
\begin{eqn}\label{eqn5.2}
\Omega_{C_1}(x)\otimes\Omega_{C_2}(x)\to
\Omega_{L_0}(P_1)\otimes\Omega_{L_1}(P_1)\otimes\cdots\otimes
\Omega_{L_{e-1}}(P_e)\otimes\Omega_{L_e}(P_e), 
\end{eqn}
sending $c(x)$ to $c(P_0)\cdots c(P_e)$. It follows that we have the
following proposition.
\begin{proposition}\label{prop5.3}
Let $x_0$ and $x_e$ be parameters of $\calO_{C_1,x}^\swedge$ and
$\calO_{C_2,x}^\swedge$, respectively. Then $c(x)$ is equal to
$dx_0\otimes dx_e$ if and only if there exist coordinates $x_i$ on the
$L_i$, for $1\leq i<e$ such that $c(P_i)=d(x_{i-1}^{-1})\otimes dx_i$
for $1\leq i\leq e$. If such coordinates exist, they are unique.
\end{proposition}
This proposition implies that knowing $X^\sim_k$ and the $c(x)$ for
the singular $x$ on it, up to isomorphism, is equivalent to knowing
$X_k$ and the $c(x)$ and $e(x)$ in that case, up to isomorphism.

The last topic of this section is the behavior of the $c(x)$ with
respect to finite morphisms. Suppose that $f\colon X\to Y$ is a finite
morphism of proper nodal curves over~$D$. Let $y$ be a singular point
of $Y_k$, and let $x$ be a point of $X_k$ that lies over~$y$. In
\cite[Prop.~5]{MaRi1} it is explained that $e(x)$ divides $e(y)$, and
that there exist isomorphisms from $D[[u,v]]/(uv-\pi^{e(x)})$ to
$\calO_{X,x}^\swedge$ and from $D[[s,t]]/(st-\pi^{e(y)})$ to
$\calO_{Y,y}^\swedge$, such that, via $f$, $s$ corresponds to $u^n$
and $t$ to $v^n$, with $n:=e(y)/e(x)$. This proves the following
proposition.
\begin{proposition}\label{prop5.4}
Let $f\colon X\to Y$ be a finite morphism of nodal curves
over~$D$. Let $y$ in $Y_k$ be a singular point, let $x$ be in
$f^{-1}y$, let $C_1$ and $C_2$ the two branches at $x$, and $D_1$ and
$D_2$ the two branches at~$y$.  Let $n$ be $e(y)/e(x)$. Then $f$ gives
isomorphisms between $\Omega_{D_i}(y)$ and~$\Omega_{C_i}(x)^{\otimes
n}$. Under these isomorphisms, $c(y)$ is mapped to~$c(x)^n$.
\end{proposition}

\section{The $c(x)$ and log structures.}\label{section6}
(This section will not be used in the rest of the article.)  The
elements $c(x)$ as described in
Sections~\ref{section3}--\ref{section5} have a simple description in
terms of log structures.  The notations are as before. Recall from
\cite{Illusie1} that a log structure on a scheme $Y$ consists of a
morphism of sheaves of multiplicative monoids $\alpha\colon
M_Y\to\calO_Y$ on the small etale site of $Y$ such that $\alpha$
induces an isomorphism from $\alpha^{-1}\calO_Y^*$ to~$\calO_Y^*$. For
our scheme $X$, let $M_X$ be the log structure induced by its special
fibre. This means that $M_X$ is the subsheaf of $\calO_X$ consisting
of functions whose restriction to the generic fibre are
invertible. Let $Q_X$ be the sheaf of monoids~$M_X/\calO_X^*$. For $x$
non-singular in $X_k(k)$, the stalk $Q_{X,x}$ is just $\NN$, with
generator the class of the uniformizer $\pi$ of~$D$. Let $x$ in
$X_k(k)$ be singular. Then, locally at $x$ in the etale topology, $X$
is isomorphic to the spectrum of $D[u,v]/(uv-\pi^{e(x)})$.  It follows
that $Q_{X,x}$ is the submonoid $(e(x),0)\NN +(0,e(x))\NN+(1,1)\NN$ of
$\NN^2$, generated by the classes of $u$, $v$ and~$\pi$. We give $D$
the log structure induced by its closed point: $M_D=D-\{0\}$. Then we
have a morphism of log schemes from $(X,M_X)$ to $(D,M_D)$. We will
now base change to $k$, equipped with its log structure $M_k$ induced
from~$M_D$: $M_k$ is the fibred sum of $M_D$ and $k^*$ over~$D^*$.
Let $(X_k,M_{X_k})$ be the pull back of $(X,M_X)$ to $(k,M_k)$. Then
$M_{X_k}$ is the fibred sum of $M_X$ and $\calO_{X_k}^*$
over~$\calO_X^*|_{X_k}$.

Let now $x$ in $X_k(k)$ be singular, and let $M_{X_k}(x)$ be the
fibred sum of $M_{X_k,x}$ and $k^*$ over~$\calO_{X_k,x}^*$. Let $C_1$
and $C_2$ be the two branches of $X_k$ at $x$, with $u=0$ on $C_1$ and
$v=0$ on~$C_2$.  We have an extension:
\begin{eqn}\label{eqn6.1}
1\to k^* \to M_{X_k}(x) \to Q_x\to 0, 
\end{eqn}
with $Q_x$ as described above. That description shows that the
$k^*$-torsors induced on the inverse images of $(e(x),0)$ and
$(0,e(x))$ correspond to the one dimensional $k$-vector spaces
$\Omega_{C_2}(x)$ and $\Omega_{C_1}(x)$, respectively. Because of the
monoid structure on $M_{X_k}(x)$, the inverse image of $(e(x),e(x))$
corresponds to $\Omega_{C_2}(x)\otimes_k\Omega_{C_1}(x)$. But via the
morphism from $(X_k,M_{X_k})$ to $(k,M_k)$ it also corresponds
to~$(m/m^2)^{\otimes e(x)}$.  The element $c(x)$ is then just the
image of~$\pi^{e(x)}$. In the case that $X$ is regular, the argument
above shows that, conversely, the $k$-scheme $X_k$ with the $c(x)$
determine the log-scheme $(X_k,M_{X_k})$ over $(k,M_k)$ up to unique
isomorphism. Compare with \cite{Steenbrink1}.

\section{The $c(x)$ for $X_0(p)$.}\label{section7}
Let $p$ be a prime number and let $X_0(p)$ and $X_0(p)^\sim$ be the
models over $\ZZ$ of modular curve $X_0(p)_\QQ$ mentioned in the
introduction.  Let $D$ be the ring of Witt vectors of an algebraic
closure $k$ of~$\FF_p$, and let $X:=X_0(p)_D$ and
$X^\sim:=X_0(p)^\sim_D$. The aim of this section is to describe the
$c(x)$ for the singular $x$ in~$X_k$. We will first give such a
description in terms of the Kodaira-Spencer map and a well known mod
$p$ modular form $B$ on supersingular elliptic curves. After that we
give an expression of some small power of $c(x)$ in terms of the
supersingular $j$-invariants mod~$p$, using the fact that the point
$0-\infty$ of $J_0(p)_\QQ$ has order $n$, the numerator of $(p-1)/12$.
Note that knowing the $c(x)$ for $X$ is equivalent to knowing them for
$X^\sim$, by Section~\ref{section5}. The results of this section will
not be used in the explicit computations of the next section; instead,
when we really need the $c(x)$, we will compute them using explicit
formulas for the Hecke correspondence~$T_2$.

Let $n\geq3$ be an integer prime to $p$, and let
$Y:=\calMbar(\Gamma(n),\Gamma_0(p))_D$.  Then $Y$ is a regular
projective nodal curve over $D$, and $Y\to X$ is the quotient for the
faithful action of the group $G:=\GL_2(\ZZ/n\ZZ)/\{1,-1\}$. The
special fibre $X_k$ of $X$ consists of two irreducible components,
$X_\infty$ and $X_0$, containing the cusps $\infty$ and $0$,
respectively. We define $Y_\infty$ to be the union of the irreducible
components of $Y_k$ that are mapped to $X_\infty$, and $Y_0$ the union
of those that are mapped to~$X_0$. Let $Z$ denote the $j$-line over
$D$, i.e., $Z:=X_0(1)_D$. The constructions that associate, to an
elliptic curve $E/S/k$, the isogenies $F\colon E\to E^{(p)}$ and
$V\colon E^{(p)}\to E$, induce isomorphisms $Z_k\to X_\infty$ and
$Z_k\to X_0$, respectively, and isomorphisms from $X(n)_k$ to
$Y_\infty$ and~$Y_0$.  The singular points of $X_k$ correspond to the
$F\colon E\to E^{(p)}$ with $E$ over $k$ supersingular. Indeed, such
an isogeny is isomorphic to $V\colon E=E^{(p^2)}\to E^{(p)}$. It
follows that the singular point $F\colon E\to E^{(p)}$ on $X_k$
corresponds to the pair $(j(E),j(E)^p)$ of points on $X_\infty$
and~$X_0$. A similar statement is of course true for~$Y_k$.

Let $P$ be a singular point of $X_k$, corresponding to a supersingular
elliptic curve $E$ over $k$, and let $Q$ be in $Y_k$, lying over~$P$.
The Kodaira-Spencer map gives an isomorphism of line bundles on
$X(n)_k$ from $\Omega(\cusps)$ to $\om^{\otimes2}$, where $\om$ is the
pullback along the zero section of the bundle of relative
differentials of the universal elliptic curve. This map, and the
isomorphisms from $Y_\infty$ and $Y_0$ to $X(n)_k$, induce
isomorphisms of $k$-vector spaces:
\begin{eqn}\label{eqn7.1}
\Omega_{Y_\infty}(Q)\to\om_E^{\otimes2},\quad
\Omega_{Y_0}(Q)\to\om_{E^{(p)}}^{\otimes2} = \om_E^{\otimes2p}. 
\end{eqn}
Via these isomorphisms, $c(Q)$ can be viewed as an element of
$\om_E^{\otimes2(p+1)}$, and then $c(P)$ is equal to~$c(Q)^{e(P)}$
in~$\om_E^{\otimes2(p+1)e(P)}$. Now there does exist, for each
supersingular elliptic curve $E$ over $k$, a non-zero element $B(E)$
in $\om_E^{\otimes(p+1)}$, with the property that for $\phi\colon E\to
E'$ any isogeny, one has $\deg(\phi)B(E)=\phi^*B(E')$. This property
determines the $B(E)$ up to a common scalar, since all supersingular
elliptic curves over $k$ are isogenous via isogenies of degree prime
to~$p$. Such $B(E)$ were defined by Robert in \cite[Thm.~B]{Robert1}
and used by Serre in \cite{Serre1}; see \cite[Prop.~7.2]{Edixhoven1}
and \cite[5]{ColemanVoloch} for several constructions. We will prove
below that, up to a common scalar, the $c(Q)$ are the squares of
the~$B(E)$. But before we do that, we give yet another construction
of~$B$.

So let $E$ be a supersingular elliptic curve over~$k$. Since $V\colon
E^{(p)}\to E$ is the dual of $F\colon E\to E^{(p)}$, we have a natural
isomorphism between the kernel $E^{(p)}[V]$ of $V$ and the Cartier
dual $E[F]^\rmD$ of the kernel of $F$ (\cite[2.8]{KatzMazur1}).  But
$E^{(p)}[V]=E^{(p)}[F]=E[F]^{(p)}$, hence we have a natural
isomorphism $\lambda_E\colon E[F]^{(p)}\to E[F]^D$. Now recall that
the category of finite commutative group schemes over $k$ that are
annihilated by $F$ and $V$ is anti-equivalent to the category of
$k$-vector spaces by sending such a group scheme to its cotangent
space at the origin, and that this converts Cartier duality into
ordinary duality (\cite[3,Thm.~7.4]{SGA}).  Hence we get an
isomorphism between $\om_E^{\otimes p}$ and $\om_E^{\otimes-1}$, i.e.,
a non-zero element $B(E)$ of~$\om_E^{\otimes(p+1)}$. For $\phi\colon
E\to E'$ an isogeny, one has
$\lambda_E\circ\phi^{(p)t}=\phi^D\circ\lambda_{E'}$, hence
$\phi^D\circ\lambda_{E'}\circ\phi^{(p)}=\deg(\phi)\lambda_E$,
establishing $\phi^*B(E')=\deg(\phi)B(E)$.
\begin{theorem}\label{theorem7.2}
There exists an element $\mu$ in $k^*$ such that for all singular $Q$
in $Y_k$ we have $c(Q)=\mu B(E)^{\otimes2}$.
\end{theorem}
\begin{proof}
Let us first note that, to define the $c(Q)$, we have to choose a
uniformizer of $D$, for which we take $p$, but that the choice does
not matter for the existence of~$\mu$. To prove the proposition, it
suffices to show that the $c(Q)$ behave in the right way with respect
to isogenies of degree prime to~$p$.

So let $Q$ be a singular point of $Y_k$, corresponding to a pair
$(E,\alpha)$, with $E$ a supersingular elliptic curve over $k$ and
$\alpha$ an isomorphism from $(\ZZ/n\ZZ)^2$ to~$E[n](k)$.  Let
$\phi\colon E\to E'$ be an isogeny of some degree $d$ that is not
divisible by~$p$. Let $Q'$ be the singular point of $Y_k$
corresponding to $(E',\phi\circ\alpha)$. Let $C$ be the modular curve
$\calM(\Gamma_0(d),\Gamma(n),\Gamma_0(p))_D$, and $s$ and $t$ the
usual morphisms from $C$ to $Y$ ($s$ just forgets the subgroup of rank
$d$, $t$ takes the quotient by it). Let $R$ be the point of $C_k$
corresponding to $(E,\phi,\alpha)$. The morphisms $s$ and $t$ are
etale, hence induce isomorphisms between $\calO_{Y,Q}^\swedge$,
$\calO_{C,R}^\swedge$ and $\calO_{Y,Q'}^\swedge$, preserving the
labeling of the branches by $\infty$ and~$0$. It follows that $c(Q)$
and $c(Q')$ correspond, via $s$ and $t$ repectively, to~$c(R)$.

For $S\to T$ be any morphism of schemes, and $f\colon F\to F'$ any
morphism of elliptic curves over $S$, the following diagram is
commutative by \cite[Lemma~7.4]{Edixhoven1}:
\begin{eqn}\label{eqn7.3}
{\renewcommand{\arraystretch}{1.5}
\begin{array}{ccc}
\om_{F/S}^{\otimes 2} & \stackrel{f^*}{\longleftarrow} & 
\om_{F'/S}^{\otimes 2} \\
{\scriptstyle {\rm KS}_{F/S}}\searrow & & 
\swarrow{\scriptstyle (\deg f){\rm KS}_{F'/S}} \\
& \Omega_{S/T} & 
\end{array}
}
\end{eqn}
Applying this to our isogeny $\phi$ at the point $R$ of $C_\infty$ and
of $C_0$ shows that the following diagram, in which the vertical
arrows are given by~(\ref{eqn7.1}), is commutative:
\begin{eqn}\label{eqn7.4}
{\renewcommand{\arraystretch}{1.5}
\begin{array}{ccc}
\Omega_{Y_\infty}(Q')\otimes\Omega_{Y_0}(Q') & 
\stackrel{(t\circ s^{-1})^*}{\longrightarrow} & 
\Omega_{Y_\infty}(Q')\otimes\Omega_{Y_0}(Q') \\
\downarrow & & \downarrow \\
\om_{E'}^{\otimes2(p+1)} & 
\stackrel{\phi^*/d^2}{\longrightarrow} & 
\om_E^{\otimes2(p+1)} 
\end{array}
}
\end{eqn}
This establishes that the $c(Q)$ behave the same as the $B(E)$ with
respect to isogenies of degree prime to~$p$.
\end{proof}
To finish this section, we will investigate the consequences of the
fact that the point $0-\infty$ in $J_0(p)(\QQ)$ has order $n$, the
numerator of the fraction $(p-1)/12$, expressed in lowest terms
(see~\cite[II, Prop.~11.1]{Mazur1}).  These consequences will not be
used in the following sections.

Let $J$ denote $\Pic^{[0]}_{X/D}/E$, as in the notation of
Section~\ref{section2}. Then $J$ is the image in the N\'eron model
$J_0(p)_D$ of~$\Pic^{[0]}_X$. Let $P$ denote the element of $J(D)$
given by $0-\infty$ (every line bundle of degree zero on $X_K$ that
admits an extension to a line bundle over $X$ determines a unique
element of~$J(D)$).  Let $e$ be the least common multiple of the
$e(x)$, $x$ ranging through the singular points of~$X_k$. We recall
from \cite[VI, Thm.~6.9]{DeligneRapoport1} that $e(x)$ equals half the
number of automorphisms of either of the two elliptic curves $E$ and
$E^{(p)}$ over $k$ corresponding to~$x$. Hence all $e(x)$ are equal to
one, except when $x$ has $j$-invariant $0$ or $1728$, in which case
$e(x)$ equals three or two, respectively, assuming $p$ is at least
five. We also recall that the elliptic curves over $k$ with
$j$-invariant equal to $0$ or $1728$ are supersingular if and only if
$p$ is congruent to $-1$ modulo three or four, respectively. It
follows that $n=(p-1)e/12$. The exact sequence (\ref{eqn2.1}) shows
that $E(D)$ is generated by the isomorphism class of the invertible
$\calO_X$-module $\calO_X(eX_\infty)$. The ``mass formula'' of
\cite[Cor.~12.4.6]{KatzMazur1} or \cite[VI, 4.9.1]{DeligneRapoport1}
gives:
\begin{eqn}\label{eqn7.5}
\deg\left(\calO_X(eX_\infty)|_{X_0}\right) = 
\sum_{x\in S_1} e/e(x) = e\sum_{x\in S_1} 1/e(x) = e(p-1)/12 = n. 
\end{eqn}
We apply Proposition~\ref{prop3.11} to the divisor $n(0-\infty)$
on~$X_k$.  It follows that there exist rational functions $f_\infty$
and $f_0$ on $X_\infty$ and $X_0$ with divisors $-n\infty+\sum_x
(e/e(x))x$ and $-n0+\sum_x (e/e(x))x$, resectively, such that at every
$x$ one has $c(x)^{e/e(x)} = f_\infty(x)f_0(x)$, where $f_\infty(x)$
and $f_0(x)$ are to be interpreted as leading terms as in
Section~\ref{section3}. This proves the following theorem.
\begin{theorem}\label{theorem7.6}
There exists a unique $\alpha$ in $\FF_p^*$ such that for every
supersingular $x$ in $X_k$ one has:
\begin{eqnarray*}
c(x)^{e/e(x)} & = & \alpha\prod_{y\neq x}(j(y)-j(x))^{e/e(y)}
\prod_{y\neq x^{(p)}}(j(y)-j(x)^p)^{e/e(y)} \cdot(dj\otimes dj)^{e/e(x)} \\
& = & \alpha\prod_{y\neq x}(j(y)-j(x))^{(p+1)e/e(y)}
\cdot(dj\otimes dj)^{e/e(x)}. 
\end{eqnarray*}
\end{theorem}
\begin{remark}\label{remark7.7}
Theorem~\ref{theorem7.6} only gives us the $c(x)^{e/e(x)}$, but
Proposition~{3.11} says that these powers are sufficient to describe
the special fibre $J_0(p)_p$ of the N\'eron model of $J_0(p)_\QQ$
completely.  Hence it should be possible to improve the main theorem
of \cite{deShalit1} in the sense that one should be able to get rid of
the factors $\pm1$ in it, because of~\cite[0.5]{deShalit1}.
\end{remark}

\section{The intersection mod $p$ of $X_0(p)$ and the cuspidal group.} 
\label{section8}
Let $p$, $k$, $X$, $X^\sim$ be as in the previous section. As in the
introduction, let $0$ and $\infty$ denote the two cusps of $X$ (we
view them as $D$-valued points of~$X$). Let $J$ denote the N\'eron
model over $D$ of $J_0(p)_K$, and let $f\colon X^{\sim,0}\to J$ be the
morphism from the smooth locus of $X^\sim$ to $J$ that extends the
morphism from $X_K$ to $J_K$ sending a point $P$ to the class of the
divisor $P-\infty$. Let $c$ be~$f(0)$. Then $c$ is of order $n$, the
numerator of $(p-1)/12$, expressed in lowest terms. Let $C$ be the
closed subgroup scheme of $J$ generated by~$c$; it is called the
cuspidal subgroup because its points correspond to classes of divisors
supported on the cusps. In this section, we want to determine
$f^{-1}C_k$, the inverse image under $f$ of the fibre over $k$ of $C$,
in order to answer the question of Coleman mentioned in the
introduction. There are two obvious elements in $f^{-1}C_k$:
$\infty_k$ and~$0_k$; the idea is that there should be nothing more.
Let us note that Mazur's theorem \cite{Mazur2} on rational isogenies
of prime degree implies that, for $p$ bigger than 163, $f^{-1}C_K$
consists just of $\infty$ and $0$, because those are the only two
$\QQ$-rational points of~$X_0(p)$. But a priori $f^{-1}C_k$ could be
bigger than the specialization of~$f^{-1}C_K$.

We will now recall some facts about the geometry of $X^\sim$, the
group $\Phi:=J_k/J_k^0$ and the map from $C_k$ to $\Phi$, that can be
found in \cite[III, Prop.~4.2]{Mazur1}.  As in the previous section,
let $X_\infty$ and $X_0$ denote the two irreducible components
of~$X_k$. The $e(x)$ for $X_k$ were described at the end of the
previous section. It follows that, for $p$ at least five, $X^\sim_k$
has at most five irreducible components: $X_\infty$ and $X_0$, one
component $F$ of genus zero (if $p$ is $-1$ mod $4$) originating from
the singularity at $x$ with $j(x)=1728$, and two components $G$ and
$H$ of genus zero (if $p$ is $-1$ mod $3$) originating from the
singularity at $x$ with $j(x)=0$. (The cases $p=2$ or $3$ will be of
no interest, because then $X_0(p)$ has genus zero.) In order to fix
the notation regarding $G$ and $H$, we demand that $G$
intersects~$X_\infty$. Figure~\ref{figure8.1} gives a picture in the
case $p=23$, where all five components are present.
\begin{figuur}\label{figure8.1} 
\centering
\caption{A picture of $X^\sim_k$ when $p=23$.}
\setlength{\unitlength}{1 cm}
\thicklines
\begin{picture}(7,8)(0,2)
\put(1,8.5){\line(1,-2){3}}
\put(6,8.5){\line(-1,-2){3}}
\put(2,6.5){\line(2,1){2}}
\put(2,6.5){\line(-2,-1){0.5}}
\put(5,6.5){\line(-2,1){2}}
\put(5,6.5){\line(2,-1){0.5}}
\put(2.5,4.5){\line(1,0){2}}
\put(1.2,8.5){$X_\infty$}
\put(6.2,8.5){$X_0$}
\put(2.7,6.5){$G$}
\put(3.9,6.5){$H$}
\put(3.3,4.6){$F$}
\end{picture}
\end{figuur} 
The intersection numbers of distinct components are as follows:
$X_\infty{\cdot}X_0$ equals the number of supersingular $j$-invariants
in $k-\{0,1728\}$; $X_\infty{\cdot}F = X_0{\cdot}F =1$;
$X_\infty{\cdot}G=1$; $X_0{\cdot}G=0$; $X_\infty{\cdot}H =0$;
$X_0{\cdot}H = 1$; $F{\cdot}G=F{\cdot}H=0$, $G{\cdot}H=1$. The self
intersections can be computed from the fact that $X^\sim_k$ has
trivial intersection with all components. From the resulting
intersection matrices one can compute, as described in
Section~\ref{section2}, that $\Phi$ is cyclic of order $n$, and
generated by the image $\phi$ of $X_0-X_\infty$ (i.e., $\Phi$ is
generated by the class of any invertible $\calO_{X^\sim_k}$-module
that has degree one on $X_0$, degree minus one on $X_\infty$, and
degree zero on the other components).  It follows that the composition
$C_k\to J_k\to\Phi$ is an isomorphism, hence every connected component
of $J_k$ contains exactly one element of~$C_k$.

Of course, if $X_K$ is of genus zero, $J$ is just the trivial group
scheme, and $f^{-1}C_k$ is all of~$X^{\sim,0}_k$. If $X_K$ is of genus
one, then $f$ is an isomorphism ($X^{\sim,0}$ and $J$ are both N\'eron
models of $X_K$), hence $f^{-1}C_k$ is isomorphic, via $f$,
to~$C_k$. So this solves the problem in these two cases, i.e., for
$p<23$. From now on we assume that $X_K$ has genus at least two, i.e.,
that $p\geq23$.
\begin{theorem}\label{theorem8.2}
Let $p\geq23$ be prime. Then we have $f^{-1}C_k=\{\infty_k,0_k\}$. 
\end{theorem}
The proof of this result will take some time. In the next section we
will see that the morphism $f\colon X^{\sim,0}_k\to J_k$ is a closed
immersion, so that $f^{-1}C_k$, as a closed subscheme, is even
reduced.  We will first prove Theorem~\ref{theorem8.2} for $p>71$, the
remaining cases will then be dealt with by hand or by simple computer
computations.

\subsection{First step of the proof of Theorem~\ref{theorem8.2}.}
\label{section8.3}
For the moment, let $p\geq23$ be prime. First we compute the map from
the set of irreducible components of $X^\sim_k$ to~$\Phi$. This
computation is done in \cite[III, Prop.~4.2]{Mazur1}, but since we
need the details of it later, we repeat it.
\begin{sublemma}\label{lemma8.3.1}
Under $f$, $X_\infty$ is mapped to $0$ in $\Phi$, $X_0$ is mapped to
$\phi$, $F$ is mapped to $(1/2)\phi$, $G$ to $(1/3)\phi$ and $H$ to
$(2/3)\phi$, (in each case the fraction makes sense mod~$n$).
\end{sublemma}
\begin{proof} Of course, $X_\infty$ 
is mapped to $0$ in $\Phi$, because $X_\infty$ contains $\infty_k$ and
$f(\infty)=0$. Likewise, $X_0$ is mapped to our generator $\phi$ of
$\Phi$, because it contains $0_k$, and $f(0)$ is the class of the
divisor $0-\infty$, which has multidegree $X_0-X_\infty$.

Suppose now that $p$ is $-1$ mod~$4$. We want to compute the image of
the component~$F$. Note that $n$ is then odd, so that we can invert
$2$ in order to compute in $\Phi$ using~(\ref{eqn2.2}). Let $P$ be
in~$F(k)$. Then $P-\infty$ has multidegree $F-X_\infty$. Note that
$\calO_{X^\sim}(F)$ has multidegree $X_\infty+X_0-2F$. By the results
of Section~\ref{section2}, the image of $f(P)$ in $\Phi$ is the class
of the multidegree
$F-X_\infty+(1/2)(X_\infty+X_0-2F)=(1/2)(X_0-X_\infty)$.

Finally, suppose that $p$ is $-1$ mod~$3$. Then $n$ is $-1$ mod~$3$.
The intersection numbers above tell us that $\calO_{X^\sim}(G)$ and
$\calO_{X^\sim}(H)$ have multidegrees $X_\infty-2G+H$ and $X_0+G-2H$,
respectively. Let $P$ be in~$G(k)$. It follows that the image of
$f(P)$ in $\Phi$ is the class of the multidegree
$$
-X_\infty+G+(1/3)\left(2(X_\infty-2G+H)+(X_0+G-2H)\right) = 
(1/3)(X_0-X_\infty). 
$$
A similar computation can be done for~$H$. 
\end{proof}
For $Z$ an irreducible component of $X_k^\sim$, let $Z^0$ be $Z$ minus
the double points of $X_k^\sim$ lying on it.
\begin{sublemma}\label{lemma8.3.2}
The only $P$ in $X_\infty^0(k)$ with $f(P)$ in $C(k)$ is~$\infty_k$.
\end{sublemma}
\begin{proof}
Suppose $P$ is in $X_\infty^0(k)$, with $f(P)$ in~$C(k)$.
Lemma~\ref{lemma8.3.1} tells us that $f(P)=0$ in~$J(k)$.  Then
Proposition~\ref{prop3.11} gives a rational function $g$ on $X_\infty$
with divisor $P-\infty$, and with a constant value at all
supersingular points. Since there are at least two distinct
supersingular points, $g$ is constant, and $P=\infty$.
\end{proof}
\begin{sublemma}\label{lemma8.3.3}
The only $P$ in $X_0^0(k)$ with $f(P)$ in $C(k)$ is~$0_k$. 
\end{sublemma}
\begin{proof}
The proof is almost the same as the previous one. In this case one
must have $f(P)=c$ in $C(k)$, and one applies
Proposition~\ref{prop3.11} to the divisor $P-0$ on~$X^\sim_k$.
\end{proof}
\begin{sublemma}\label{lemma8.3.4}
Let $p\geq23$ be $-1$ mod $4$ and suppose that there exists a
supersingular $j$-invariant in $\FF_p-\{1728\}$ and one in
$\FF_{p^2}-\FF_p$. Then there is no point $P$ in $F^0(k)$ such that
$2f(P)$ is in~$C(k)$.
\end{sublemma}
\begin{proof}
Suppose that $p$ is as in the Lemma, and that $P$ in $F^0(k)$ has
$2f(P)$ in~$C(k)$. Then $2f(P)=c$ by Lemma~\ref{lemma8.3.1}.  We apply
Proposition~\ref{prop3.11} to the divisor $E:=2P-0-\infty$. It follows
that there exist $\lambda_\infty$ and $\lambda_0$ in $k^*$, such that
the functions $f_\infty=\lambda_\infty(j-1728)$ and
$f_0=\lambda_0(j-1728)$ on $X_\infty$ and $X_0$ satisfy
$f_\infty(x)=f_0(x^p)$ for every supersingular $j$-invariant $x$ in
$k$ other than~$1728$. The existence of such an $x$ in $\FF_p$ shows
that $\lambda_\infty=\lambda_0$. But the existence of such an $x$ in
$\FF_{p^2}-\FF_p$ implies that $\lambda_\infty=-\lambda_0$.
\end{proof}
\begin{sublemma}\label{lemma8.3.5}
Let $p\geq23$ be $-1$ mod $3$ and suppose that there exist two
supersingular $j$-invariants $x$ and $y$ in $k-\{0\}$ such that
$x^{p-2}\neq y^{p-2}$.  Then there is no point $P$ in $G^0(k)$ or in
$H^0(k)$ such that $3f(P)$ is in~$C(k)$.
\end{sublemma}
\begin{proof}
Suppose that $p$ is as in the Lemma, and that $P$ in $G^0(k)$ has
$3f(P)$ in~$C(k)$. Then $3f(P)=c$ by Lemma~\ref{lemma8.3.1}. We apply
Proposition~\ref{prop3.11} to the divisor $E:=3P-2\infty-0$. It
follows that $\calO_{X^\sim}(-2G-H)|_{X^\sim_k}$ has a rational
section with divisor~$E$. The restrictions to $X_\infty$ and $X_0$ of
this section are of the form $\lambda_\infty j^2$ and $\lambda_0j$,
respectively. It follows that for every supersingular $x$ in $k-\{0\}$
we have $\lambda_\infty x^2=\lambda_0 x^p$. This contradicts the
hypotheses of the Lemma. The same result holds for $H$ because of the
$w_p$-operator.
\end{proof}
\begin{sublemma}\label{lemma8.3.6}
For every prime $p\geq23$ with $p=-1$ mod $3$ there do exist
supersingular $j$-invariants $x$ and $y$ in $k-\{0\}$ with
$x^{p-2}\neq y^{p-2}$.
\end{sublemma}
\begin{proof}
Note that $p-2$ and $p^2-1$ have greatest common divisor $3$, hence
the required $x$ and $y$ exist if there are at least four distinct
non-zero supersingular $j$-invariants in~$k$. The mass formula of
\cite[Cor.~12.4.6]{KatzMazur1} shows that the number of supersingular
$j$-invariants is at least $(p-1)/12$, which proves the Lemma for
$p>49$.  The primes that are left are $23$, $29$, $41$
and~$47$. Table~6 in \cite{MFIV} shows that in each of these cases
there are at least two non-zero supersingular $j$-invariants $x$ and
$y$ in~$\FF_p$. For these we have $x^{p-2}=x^{-1}\neq y^{-1}=y^{p-2}$.
\end{proof}
\begin{sublemma}\label{lemma8.3.7}
For every prime $p>7$ with $p=-1$ mod $4$ there exists a supersingular
$j$-invariant in $\FF_p-\{1728\}$.
\end{sublemma}
\begin{proof}
There is a unique elliptic curve over $\CC$ (up to isomorphism) with
endomorphism ring isomorphic to~$\ZZ[2i]$. Its $j$-invariant is
$2^33^311^3$, and we have $2^33^311^3-1728=2^33^67^2$. (An easy way to
compute this is to use the explicit description of the Hecke
correspondence $T_2$ given below in Section~\ref{section8.4}.)
\end{proof}
\begin{sublemma}\label{lemma8.3.8}
For every prime $p>71$ there exists a supersingular $j$-invariant in
$k-\FF_p$.
\end{sublemma}
\begin{proof}
This is a result of Ogg, which can be found as a combination of the
identities~(21) and (24) in~\cite{Ogg1}. We recall the argument in a
few lines. The $\FF_p$-rational supersingular points of $X_k$ are
precisely the fixed points of the Fricke involution $w$ of~$X$. It
follows that all supersingular $j$-invariants in $k$ are in $\FF_p$ if
and only if $X_k^+:=X_k/\{1,w\}$ has genus zero, hence if and only if
$X_0(p)/\{1,w\}$ is isomorphic to~$\PP^1_\ZZ$. Suppose that this is
the case. Then the inequalities:
$$
10 = 2\#\PP^1(\FF_4)\geq \#X_0(p)(\FF_4)\geq (p+1)/12 + 2
$$
show that $p$ is at most~$95$. Table~6 of \cite{MFIV} then finishes
the proof.
\end{proof}
Combining the previous lemmas in this section proves
Theorem~\ref{theorem8.2} for all primes $p>71$:
Lemmas~\ref{lemma8.3.5} and \ref{lemma8.3.6} take care of points on
$G$ and $H$, and Lemmas~\ref{lemma8.3.4}, \ref{lemma8.3.7} and
\ref{lemma8.3.8} take care of points on~$F$.  For the remaining primes
$p$ with $23\leq p\leq71$ with $p=-1$ mod $4$ we still have to show
that no point $P$ in $F^0(k)$ has $f(P)$ in~$C(k)$. Another inspection
of Table~6 of \cite{MFIV} shows that the hypothesis of
Lemma~\ref{lemma8.3.4} is satisfied for $43$ and~$67$.  The only
primes with which we still have to deal are $23$, $31$, $47$, $59$
and~$71$.

\subsection{Second step in the proof of Theorem~\ref{theorem8.2}.}
\label{section8.4}
In this section we will deal with the remaining $p$ for which we still
have to prove Theorem~\ref{theorem8.2}: $23$, $31$, $47$, $59$
and~$71$. We recall that in these cases we only have to show that no
$P$ in $F^0(k)$ has $f(P)$ in~$C(k)$.  The proof will use the full
strength of Proposition~\ref{prop3.11} and the behavior of the $c(x)$
under the Hecke correspondence~$T_2$. Until now we only applied
Proposition~\ref{prop3.11} in rather simple cases, not involving
the~$c(x)$.

In order to apply Proposition~\ref{prop3.11}, it suffices to know the
$c(x)$ up to a common scalar, for example, since the choice of a
uniformizer of $D$ is irrelevant. To compute the $c(x)$ up to a common
scalar, we just pick a supersingular point $x_0$ with $e(x_0)=1$ and
choose $c(x_0)=dj\otimes dj$. Then we use the Hecke correspondence
$T_2$ on $X$ to find the $c(x)$ at the other supersingular~$x$. The
following explicit description of the Hecke correspondence $T_2$ on
the affine $j$-line (over $\ZZ$) is given in~\cite{Mestre1}:
$Y_0(2)_\ZZ$ is the curve in $\AA^2_\ZZ$ given by the equation
$uv=2^{12}$, the morphisms $s$ and $t$ from $Y_0(2)_\ZZ$ to the
$j$-line are given by $s^*j=(u+16)^3/u$ and $t^*j=(v+16)^3/v$. Now
consider the correspondence $T_2$ on $X$; it is given by the morphisms
$s$ and $t$ from $X_0(2p)_D$ to~$X$. In the beginning of
Section~\ref{section7} we gave isomorphisms between the $j$-line
$\PP^1_k$ and the two irreducible components $X_\infty$ and $X_0$
of~$X_k$. Likewise, $X_0(2p)_k$ is the union of $X_0(2p)_\infty$ and
$X_0(2p)_0$, both isomorphic to $X_0(2)_k$ via similar isomorphisms,
and the induced correspondences on $X_\infty$ and $X_0$ are given by
the explicit formulas above.
\begin{sublemma}\label{lemma8.4.1}
Let $y$ in $X_0(2p)(k)$ be supersingular, and put $x:=s(y)$,
$x':=t(y)$.  Then $j-j(x)$ is a local coordinate on $X_\infty$ at $x$,
$j-j(x)^p$ is one on $X_0$ at $x$, $u-u(y)$ and $v-v(y)$ are
coordinates on $X_0(2p)_\infty$ at $y$, $u-u(y)^p$ and $v-v(y)^p$ are
coordinates on $X_0(2p)_0$ at $y$, and $j-j(x')$ and $j-j(x')^p$ are
local coordinates at $x'$ on $X_\infty$ and $X_0$,
respectively. Recall from Proposition~\ref{prop5.4} that $e(y)$
divides $e(x)$ and~$e(x')$. Let $\alpha$ in $k^*$ be such that
$c(y)=\alpha (du\otimes du)$. Let $\beta$ and $\gamma$ in $k^*$ be
given by:
\begin{eqnarray*}
s^*(j-j(x)) & = & \beta (u-u(y))^{e(x)/e(y)} + \mbox{higher order terms,}\\
t^*(j-j(x')) & = & \gamma (v-v(y))^{e(x')/e(y)} + \mbox{higher order terms.}
\end{eqnarray*}
Then we have: 
\begin{eqnarray*}
c(x) & = & \alpha^{e(x)/e(y)}\beta^{-(p+1)}dj\otimes dj, \\
c(x') & = & (-\alpha u(y)/v(y))^{e(x')/e(y)}\gamma^{-(p+1)}dj\otimes dj. 
\end{eqnarray*}
\end{sublemma}
\begin{proof} 
The local coordinates are obtained using the isomorphisms from the
$j$-line to $X_\infty$ and $X_0$ and from $X_0(2)_k$ to
$X_0(2p)_\infty$ and $X_0(2p)_0$ that we discussed above.  The rest is
a direct application of Proposition~\ref{prop5.4}.
\end{proof}
\begin{subremark}\label{remark8.4.2}
In practice, we will start with a supersingular $x$ in $X(k)$ with
$e(x)=1$.  In order to find the various $y$'s in $s^{-1}x$, we factor
the polynomial $(T+16)^3-j(x)T$ modulo~$p$. (Actually, we will always
have $j(x)$ in~$\FF_p$.) The factor $\beta$ of Lemma~\ref{lemma8.4.1}
is then just the derivative $(3(u(y)+16)^2-j(x))/u(y)$. Let us write
$c(x)=\delta\, dj\otimes dj$. Suppose that $e(x)=1=e(x')$. Then we
simply have:
$$
c(x')=\left(\frac{j(x')-3(v(y)+16)^2}{3(u(y)+16)^2-j(x)}\right)^{p+1}
\delta\, dj\otimes dj.
$$
In all our computations, $j(x)$, $j(x')$, $u(y)$ and $v(y)$ will be in
$\FF_p$, so that we can replace the exponent $p+1$ by~$2$.
\end{subremark}

\subsubsection{The case $p=31$.}\label{section8.4.3} 
The supersingular $j$-values in $k$ are $1728=-8$, $2$ and $4$ (see
for example Table~6 in~\cite{MFIV}).  We have to check that no
non-singular point on the component $F$ is mapped to the cuspidal
group. Suppose that $P$ is in $F^0(k)$ with $f(P)$ in~$C(k)$.
Lemma~\ref{lemma8.3.1} shows that $f(P)=3c$. Hence we apply
Proposition~\ref{prop3.11} to the divisor $-2\infty+3{\cdot}0-P$. Note
that $\calO_{X^\sim}(X_0)$ has multidegree $2X_\infty+F-2X_0$. Let
$x_2$ and $x_4$ denote the supersingular points of $X(k)$ with
$j$-invariants $2$ and $4$, respectively. It follows that the
functions $f_\infty=(j-2)(j-4)$ on $X_\infty$, $f_0=1/(j+8)(j-2)(j-4)$
on $X_0$ satisfy equation~\ref{eqn3.12} for the cycle
$(X_\infty,x_2,X_0,x_4,X_\infty)$.  As we remarked above, we only need
to know $c(x_2)$ and $c(x_4)$ up to a common scalar, hence we simply
put $c(x_2)=dj\otimes dj$. We apply the method described in
Remark~\ref{remark8.4.2}. The $y$'s in $s^{-1}x_2$ are given by
$u(y)\in\{-15,2,-4\}$. We take $y$ with $u(y)=-4$; then
$j(t(y))=4$. It follows that $c(x_4)=-6dj\otimes dj$. The factor at
$x_2$ in equation~\ref{eqn3.12} is:
$$
((2-4)dj\otimes((2+8)(2-4))dj)c(x_2)^{-1}=9(dj\otimes dj)c(x_2)^{-1}=9. 
$$
The analogous factor at $x_4$ is: 
$$
((4+8)(4-2)dj\otimes(4-2)dj)^{-1}c(x_4)=
-6{\cdot}17^{-1}(dj\otimes dj)^{-1}(dj\otimes dj)=-4. 
$$
Since $-4{\cdot}9=-5\neq1$ in $\FF_{31}$, we have found a
contradiction, so we have proved that no $P$ in $F^0(k)$ has image
in~$C(k)$.

\subsubsection{The case $p=47$.}\label{section8.4.4}
The supersingular $j$-invariants are $0$, $1728=-11$, $9$, $10$
and~$-3$.  Suppose that $P$ is in $F^0(k)$ with $f(P)$ in~$C(k)$. Then
$f(P)=12c$ by Lemma~\ref{lemma8.3.1}. We apply
Proposition~\ref{prop3.11} to the divisor
$E:=-11\infty+12{\cdot}0-P$. The divisor $X_\infty-2X_0-H$ has the
same multidegree as~$E$. Let $x_9$ and $x_{10}$ in $X(k)$ be the
points with $j$-invariants $9$ and~$10$. The functions:
$$
\mbox{$f_\infty=j(j+11)(j-9)^3(j-10)^3(j+3)^3$ and 
$f_0=(j(j+11)^2(j-9)^3(j-10)^3(j+3)^3)^{-1}$}
$$ 
then satisfy equation~\ref{eqn3.12} for the cycle
$(X_\infty,x_9,X_0,x_{10},X_\infty)$. We normalize the $c(x)$ such
that $c(x_9)=dj\otimes dj$. In order to compute $c(x_{10})$, consider
the point $y$ on $X_0(2)(k)$ with $u(y)=-1$. Then we have $j(s(y))=9$
and $j(t(y))=10$. One computes that $c(x_{10})=3dj\otimes dj$, and
that the factors at $x_9$ and $x_{10}$ in equation~\ref{eqn3.12} equal
$15$ and $-11$, respectively. Since $-11{\cdot}15=23\neq1$, we have
proved that no $P$ in $F^0(k)$ has $f(P)$ in~$C(k)$.

\subsubsection{The case $p=59$.}\label{section8.4.5}
The supersingular $j$-invariants are $0$, $1728=17$, $15$, $28$, $-12$
and~$-11$.  Suppose that $P$ is in $F^0(k)$ with $f(P)$
in~$C(k)$. Then computations as in the previous two cases show that
the functions:
\begin{eqnarray*}
f_\infty & = & j(j-17)(j-15)^3(j-28)^3(j+12)^3(j+11)^3 \\
f_0 & = & (j(j-17)^2(j-15)^3(j-28)^3(j+12)^3(j+11)^3)^{-1} 
\end{eqnarray*} 
satisfy equation~\ref{eqn3.12} for the cycle
$(X_\infty,x_{15},X_0,x_{28},X_\infty)$. We normalize the $c(x)$ such
that $c(x_{15})=dj\otimes dj$. In order to compute $c(x_{28})$,
consider the point $y$ on $X_0(2)(k)$ with $u(y)=-3$. One finds that
$c(x_{28})=21dj\otimes dj$, and that the factors at $x_{15}$ and
$x_{28}$ in equation~\ref{eqn3.12} equal $35$ and $10$,
respectively. Since $35{\cdot}10=-4\neq1$, we have proved that no $P$
in $F^0(k)$ has $f(P)$ in~$C(k)$.

\subsubsection{The case $p=71$.}\label{section8.4.6}
The supersingular $j$-invariants are $0$, $1728=24$, $17$, $-31$,
$-30$, $-23$ and~$-5$.  Suppose that $P$ is in $F^0(k)$ with $f(P)$
in~$C(k)$. Computations as in the previous cases show that the
functions:
\begin{eqnarray*}
f_\infty & = & j(j-24)(j-17)^3(j+31)^3(j+30)^3(j+23)^3(j+5)^3 \\
f_0 & = & (j(j-24)^2(j-17)^3(j+31)^3(j+30)^3(j+23)^3(j+5)^3)^{-1} 
\end{eqnarray*} 
satisfy equation~\ref{eqn3.12} for the cycle
$(X_\infty,x_{-31},X_0,x_{-23},X_\infty)$. We normalize the $c(x)$
such that $c(x_{-31})=dj\otimes dj$, and find that
$c(x_{-23})=-13dj\otimes dj$ using the point $y$ with $u(y)=-14$.  The
factors at $x_{-31}$ and $x_{-23}$ in equation~\ref{eqn3.12} equal
$49$ and $38$, respectively. Since $49{\cdot}38=16\neq1$, we have
proved that no $P$ in $F^0(k)$ has $f(P)$ in~$C(k)$.

\subsubsection{The case $p=23$.}\label{section8.4.7}
The supersingular $j$-invariants are $0$, $1728=3$ and~$-4$.  Suppose
that $P$ is in $F^0(k)$ with $f(P)$ in~$C(k)$. In this case we work
with the model $X'$ of $X$ obtained by blowing up $X$ in its
$k$-valued point with $j$-invariant~$1728$. We apply
Proposition~\ref{prop3.11} to the divisor $-5\infty+6{\cdot}0-P$. For
$a$ in the proposition one can take $-3X_0-F$. One finds that the
functions:
$$
\mbox{$f_\infty=j(j-3)(j+4)^3$ and $f_0=(j(j-3)^2(j+4)^3)^{-1}$} 
$$
satisfy equation~\ref{eqn3.12} for the cycle
$(X_\infty,x_0,X_0,x_{-4},X_\infty)$. We normalize the $c(x)$ such
that $c(x_{-4})=dj\otimes dj$. Consider the point $y$ with
$u(y)=-3$. In the notation of Lemma~\ref{lemma8.4.1} one finds:
$\beta=6$, $\alpha=\beta^2=13$ and $\gamma=1/7$. It follows that
$c(x_0)=-3dj\otimes dj$.  The factors at $x_0$ and $x_{-4}$ in
equation~\ref{eqn3.12} equal $13$ and $-5$, respectively. Since
$13{\cdot}-5=4\neq1$, we have proved that no $P$ in $F^0(k)$ has
$f(P)$ in~$C(k)$.

\section{Closed immersions.} \label{section9}
We go back to the situation of Section~\ref{section2}. So let $D$ be a
complete discrete valuation ring with fraction field $K$, separably
closed residue field $k$, and uniformizer~$\pi$. Let $X_K$ be a proper
smooth geometrically irreducible curve, and $X$ a proper flat nodal
model of it.  Let $J:=\Pic^{[0]}_{X/D}/E$ be the open part of the
N\'eron model $J^\sim$ of the jacobian of $X_K$ over $D$ corresponding
to line bundles on $X_K$ that admit an extension to a line bundle
on~$X$. Let $X^0$ be the open part of $X$ where the morphism to
$\Spec(D)$ is smooth, i.e., $X^0$ is the complement in $X$ of the set
of singular points of~$X_k$. Let $P$ be in $X^0(D)$. Then we have a
morphism $f\colon X^0\to J$ that sends, for $S$ any $D$-scheme, $Q$ in
$X^0(S)$ to the class of $\calO_{X_S}(Q-P_S)$. Equivalently, $f$ is
the morphism one obtains by applying the N\'eron property to $X^0$ and
$f_K\colon X_K\to J_K$, using that the image of $X^0$ in $J^\sim$ is
in~$J$. In this section we study to what extent $f$ is a closed
immersion. The first observation in doing this is that the
disconnecting double points behave differently from the
non-disconnecting ones. Recall that a singular point $x$ of $X_k$ is
called disconnecting if and only if $X_k-\{x\}$ is not connected. We
let $X'$ be the complement in $X$ of the set of non-disconnecting
double points of~$X_k$.
\begin{proposition}\label{prop9.1}
The morphism $f\colon X^0\to J$ extends to a morphism $f\colon X'\to J$. 
\end{proposition}
\begin{proof}
Let $i\colon\Gamma\to X^0\times_DJ$ be the graph of $f$, and let
$\ol{\Gamma}$ be its closure in $X'\times_DJ$. We have to show that
the projection $p$ from $\ol{\Gamma}$ to $X'$ is an isomorphism. We
will show, using the valuative criterion of properness (\cite[II,
7.3.8 and 7.3.9]{EGA}), that $p$ is proper. Then we show that $p$ is
bijective, hence finite. The normality of $X$ then implies that $p$ is
an isomorphism.

In order to prove that $p$ is proper, it suffices to show that for
$D'$ any complete discrete valuation ring with fraction field $K'$ and
algebraically closed residue field $k'$, and any $Q$ in $X'(D')$ with
$Q(\Spec(K'))$ in $X_K$, the element $f(Q)$ in $J(K')$ extends
uniquely to one in~$J(D')$.  (The fact that it is sufficient to
consider $Q$ in $X'(D')$ with $Q(\Spec(K'))$ in the open dense subset
$X_K$ of $X'$ can be found for example, without proof, in
\cite[p.~103]{DeligneMumford}, and, with a proof, in the forthcoming
book~\cite{LaumonMoretBailly}.  The proof of this fact starts with a
reduction to the quasi-projective case, using Chow's Lemma; then one
considers a compactification.)  Let $D'$ and $Q$ in $X'(D')$ be as
above. If $Q$ is in $X^0(D')$, the condition is obviously satisfied
(just compose $Q$ with~$f$). So we suppose that $Q$ is not
in~$X^0(D')$. That means that $Q(\Spec(k'))$ is a disconnecting double
point, say $x$, of~$X_k$. Base change from $D$ to $D'$ reduces the
whole problem to the case $D'=D$. By definition, $f_K(Q)$ extends in
$J^\sim$; we have to show that that extension is actually
in~$J$. Equivalently, we have to show that the connected component of
$J^\sim_k$ to which $f_K(Q)$ specializes is a connected component
of~$J_k$. Let $F_1,\ldots,F_{e(x)-1}$ be the chain of $\PP^1_k$'s
lying over $x$ in the minimal resolution $R\colon X^\sim\to X$. Since
$x$ is disconnecting, $X^\sim_k$ can be written as $Y\cup
Z\cup(\cup_iF_i)$, with $Y$ and $Z$ connected, disjoint, and $Y\cup Z$
the union of all irreducible components of $X^\sim_k$ other than
the~$F_i$. To fix who is who, we demand that $P$ is in $Y$ and that
$F_1$ meets~$Y$. In $X^\sim$, $Q$ specializes to one of the $F_i$, say
$F_i$. Let $C$ be the irreducible component of $X_k$ to which $P$
specializes. Then the connected component of $J^\sim_k$ to which
$f(Q)$ specializes is, in the description of Section~\ref{section2},
the class of the multidegree $F_i-C$. Let $a$ be the divisor on
$X^\sim$ which is given as the sum of all irreducible components of
$Y$, and the $F_j$ with $j<i$. The multidegree of $a$ equals
$F_i-F_{i-1}$. It follows that the class of the multidegree $F_i-C$ is
the same as that of $F_0-C$, if $F_0$ denotes the unique irreducible
component of $Y$ meeting~$F_1$. Since $F_0$ maps onto an irreducible
component of $X_k$, we have shown that the class of $F_i-C$ is
represented by a multidegree with support in the set of irreducible
components of $X_k$, hence that $f(Q)$ extends to an element
of~$J(D)$. This finishes the proof that $p$ is proper.

Let us now show that $p$ is bijective. So let $x$ be a disconnecting
double point of~$X_k$. It suffices now to show that for every $Q$ in
$X(K)$ specializing to $x$, $f(Q)$ specializes to the same element
in~$J(k)$ (the argument is similar to the proof of the generalization
in \cite{DeligneMumford} of the valuative criterion of properness
mentioned above). Recall that $J_k=\Pic^{[0]}_{X_k/k}/E_k$. Let $Q$
and $Q'$ be two elements of $X(K)$ specializing to~$x$. The fact that
the $F_i$ are isomorphic to $\PP^1_k$ implies that $f(Q)$ and $f(Q')$
specialize to the same element of $J(k)$ if $Q$ and $Q'$ specialize to
the same $F_i$; one reduces to that situation using divisors of the
type $Y+\sum_{j<i}F_j$ as above.
\end{proof}
Now that we know that $f$ extends to $f\colon X'\to J$, we want to
know when this extension is a closed immersion. We begin by analyzing
what happens in the special fibres. We denote by $G$ the graph of
$X_k$ as in Section~\ref{section2}, with $S_0$ as set of vertices and
$S_1$ as set of edges. We use the analogous notation $G^\sim$,
$S_0^\sim$ and $S_1^\sim$ for~$X^\sim_k$.
\begin{proposition}\label{prop9.2}
Suppose that $C_1$ and $C_2$ are two distinct elements of $S_0$ such
that the multidegree $C_1-C_2$ has image zero in
$\Phi=J_k/J_k^0$. Then there is a unique path in $G$ from $C_1$ to
$C_2$ such that every edge in that path is a disconnecting double
point.
\end{proposition}
\begin{proof}
The graph $G^\sim$ is obtained from $G$ by replacing each edge $x$ by
a chain of $e(x)$ edges, hence it suffices to prove this proposition
in the case $X=X^\sim$. So we assume that $X$ is regular. We choose an
orientation on $G$, i.e., maps $s$ and $t$ from $S_1$ to $S_0$ such
that, for every $x$ in $S_1$, $s(x)$ and $t(x)$ are the vertices
of~$x$. Then we have maps $s_*$ and $t_*$ from $\ZZ^{S_1}$ to
$\ZZ^{S_0}$ and $s^*$ and $t^*$ from $\ZZ^{S_0}$ to $\ZZ^{S_1}$ given
by $(s_*f)(C)=\sum_{s(x)=C}f(x)$ and $(s^*f)x=f(s(x))$, etc. We define
the usual boundary and coboundary maps $d_*:=t_*-s_*$ and
$d^*:=t^*-s^*$.  We recall from \cite[1]{Edixhoven2} that we have the
exact sequence:
\begin{eqn}\label{eqn9.3}
0\lto \ZZ\;\;\stackrel{\rm diag}{\lto}\;\;
\ZZ^{S_0}\;\;\stackrel{d_*d^*}{\lto}\;\;\ZZ^{S_0}[+] \lto \Phi \lto 0, 
\end{eqn}
with $\ZZ^{S_0}[+]$ the kernel of $\ZZ^{S_0}\to\ZZ$,
$f\mapsto\sum_Cf(C)$.  In fact, this exact sequence is just a simple
reformulation of Raynaud's description of $\Phi$ from
Section~\ref{section2}. (The hypothesis in \cite{Edixhoven2} that the
irreducible components of $X_k$ are smooth is actually never used in
that article, and the hypothesis that $k$ be algebraically closed can
be replaced with $k$ separably closed.)  All this is to indicate that
$\Phi$ is the cokernel of a kind of Laplace operator. For $g$ in
$\ZZ^{S_0}[+]$, there exists $h$ in $\QQ^{S_0}$, unique up to adding a
constant function, such that $g=d_*d^*h$. Such functions $g$ and $h$
have the following electrodynamical interpretation.  We view $G$ as an
electric circuit in which each edge has a resistance of one ohm, and
in which at each $C$ in $S_0$ an electric current of $g(C)$ amp\`ere
enters the circuit; then $h$ is a function that gives at each $C$ the
potential in volts.

In our situation, we know that for $g=C_1-C_2$ there exists such an
$h$ with integer values. Since the total current leaving the circuit
is one, the current in each edge is then $0$, $1$ or~$-1$. This
implies that if we delete an edge through which there goes a non-zero
current, the resulting graph is disconnected, which gives the
conclusion we wanted.

For those readers (including the referee) who do not accept such
physical arguments as mathematics, we give the following rigorous
argument.  Let $V$ and $E$ be the sets of vertices and edges. Say one
amp\`ere enters at $x$ and leaves at $y$.  Let $V_{\rm max}$ be the
set of $z$ in $V$ such that $h(z)=h(x)$ and $z\neq x$. Let $V_{\rm
min}$ be the set of $z$ in $V$ such that $h(z)=h(y)$ and $z\neq
y$. Then $V-\{x,y\}$ is the disjoint union of $V_{\rm max}$, $V_{\rm
min}$ and say~$V'$.  Since $d_*d^*h$ is zero at all $z$ in
$V-\{x,y\}$, it follows that $V_{\rm max}$, $V_{\rm min}$ and $V'$ are
disconnected from each other in~$G-\{x,y\}$.  Since every edge
connecting $x$ to $V'$ carries a non-zero current, there is a unique
such edge. Likewise for~$y$. Induction on the number of elements of
$V$ (remove $V_{\rm max}$ and $V_{\rm min}$, and $x$ and~$y$) finishes
the proof.
\end{proof}
\begin{proposition}\label{prop9.4}
The morphism $f_k\colon X'_k\to J_k$ is proper. 
\end{proposition}
\begin{proof}
We apply the valuative criterion of properness. It suffices to
consider the local rings of $X_k^\nor$ at closed points $x'$ whose
image $x$ in $X_k$ is a non-disconnecting double point. Let $x$ be
such a point. We have to show that $f_k$ does not extend over some
neighborhood of~$x$. Suppose first that $x$ is an intersection point
of two distinct irreducible components $C$ and $C'$
of~$X_k$. Proposition~\ref{prop9.2} implies that the images of $C\cap
X'_k$ and $C'\cap X'_k$ under $f_k$ lie in two distinct connected
components of~$J_k$. It follows that $f_k$ does not extend
over~$x$. To finish, suppose that $x$ is a double point of an
irreducible component $C$ of~$X_k$. Let $C'$ be the partial
normalization of $C$, in which all singular points are normalized
except~$x$. As before, we may and do suppose that $P$ specializes
to~$C$. Then $f_k(C\cap X'_k)$ lies in $J_k^0=\Pic^0_{X_k/k}$.  Let
$g\colon J_k^0\to\Pic^0_{C'/k}$ be the morphism obtained by
restriction of line bundles. We claim that $g\circ f_k$ does not
extend over~$x$.  Suppose that it does. Then $\Pic^0_{C'/k}$ would be
generated by the image of $C'$ in it, hence would be a complete
variety, which it isn't (it is an extension of $\Pic^0_{C^\nor/k}$
by~$\GG_\rmm$).
\end{proof}

\begin{proposition}\label{prop9.5}
The following two conditions are equivalent: 
\begin{enumerate}
\item $f_k\colon X'_k\to J_k$ is injective; 
\item every irreducible component of $X_k$ isomorphic to $\PP^1_k$ 
contains a non-disconnecting double point. 
\end{enumerate}
If these conditions are satisfied, then $f_k$ is a closed immersion. 
\end{proposition}
\begin{proof}
Let us first show that the first condition implies the second. So
assume that the second condition is not satisfied, and let $C$ is an
irreducible component of $X_k$, isomorphic to $\PP^1_k$ and not
containing any non-disconnecting double point. Let $C^0$ be denote the
complement in $C$ of the set of double points of $X_k$ contained
in~$C$. Let $Q$ and $Q'$ be two distinct elements of~$C^0(k)$. An easy
application of Proposition~\ref{prop3.11} shows that $f(Q)=f(Q')$.

Assume now that the second condition is satisfied. We show that $f_k$
is a closed immersion. Since we know that $f_k$ is proper, it suffices
to show that $f_k$ is injective, and injective on tangent spaces. Let
us first consider the restrictions of $f_k$ to the irreducible
components of $X'_k$, and prove that these are injective and injective
on tangent spaces.

Let $C$ be an irreducible component of~$X_k$. Replacing $P$ by an
element of $X^0(D)$ that specializes to $C$ changes $f$ by a
translation, hence we may and do assume that $P$ specializes
to~$C$. Then $f_k(C\cap X'_k)$ is contained in
$J^0_k=\Pic^0_{X_k/k}$. Since the double points of $X_k$ contained in
$C\cap X'_k$ are disconnecting they are smooth points of $C$, hence
the morphism $f_k\colon C\cap X'_k\to J^0_k$ is the unique extension
of its restriction to~$C^0$. If $C$ is not isomorphic to $\PP^1_k$
then the composition of $f_k\colon C\cap X'_k\to J^0_k$ with the
morphism $J^0_k\to\Pic^0_{C/k}$ induced by restriction of line bundles
on $X_k$ to $C$ is injective and is injective on tangent spaces, hence
the same is true for $f_k\colon C\cap X'_k\to J^0_k$.
So suppose that $C$ is isomorphic to~$\PP^1_k$. Then $f_k(C\cap X'_k)$
is contained in the maximal torus $T$
of~$J_k^0$. Since $C$ does contain a non-disconnecting double point,
there exists a connected component $D$ of $X_k-C$ whose closure in
$X_k$ meets $C$ in at least two points, say $x$ and~$y$. Let $g$ be
the regular function on $C-\{x,y\}$ with divisor $x-y$ and
$g(P_k)=1$. The standard description of $T$ says that the character
group of $T$ is the homology group $\rH_1(G,\ZZ)$ of the graph~$G$.
It follows that there is a character $\chi\colon T\to\GG_\rmm$ such
that $f_k\colon C\cap X'_k\to T$ composed with $\chi$ is the
function~$g$. Hence also in this case the restriction of $f$ is
injective and injective on tangent spaces.

To finish the proof, it suffices to prove that $f_k$ is injective, and
injective on the tangent spaces at the double points
of~$X_k'$. Suppose that $Q$ and $Q'$ in $X'(k)$ have
$f_k(Q)=f_k(Q')$. Proposition~\ref{prop9.2} then implies that $Q$ and
$Q'$ lie on irreducible components $C$ and $C'$ of $X_k$ that are
connected in $G$ by a path consisting of disconnecting double
points. Let $C_1,\ldots,C_n$ be the irreducible components in this
path, and let $Y$ denote their union, as a reduced closed subscheme
of~$X_k$.  After replacing $P$ by a point that specializes to $Y\cap
X^0$, $f_k$ sends $Y\cap X'$ into~$J^0_k$. Let $h\colon
J^0_k\to\Pic^0_{Y/k}$ be the morphism obtained by restriction of line
bundles. Since the double points in the path above are disconnecting,
the morphism from $\Pic^0_{Y/k}$ to the product of the
$\Pic^0_{C_i/k}$ is an isomorphism. For each $i$, let $p_i$ be the
projection from $\Pic^0_{Y/k}$ to~$\Pic^0_{C_i/k}$.
Proposition~\ref{prop3.11} implies that the restriction of $p_ihf_k$
to $C_j\cap X'$ is constant if $i\neq j$, and an immersion if
$i=j$. This proves that $Q=Q'$, and also that $f_k$ is injective on
the tangent spaces at the double points.
\end{proof}
\begin{proposition}\label{prop9.6}
The following three conditions are equivalent: 
\begin{enumerate}
\item the morphism $f\colon X'\to J$ is proper; 
\item for every prime number $p$, every non-disconnecting double point
$x$ of $X_k$ is contained in a loop in the graph $G$ such that
$v_p(e(x))\leq v_p(e(y))$ for every $y$ in that loop, where $v_p$
denotes the $p$-adic valuation;
\item for every finite extension $K'$ of $K$ and every $Q$ of $X(K')$
specializing to a non-dis\-connecting double point, $f(Q)$ specializes
to a connected component of $J_{D'}^\sim$ that is not in $J_{D'}$,
with $D'$ denoting the integral closure of $D$ in $K'$, $J_{D'}$ the
pullback of $J$ to $D'$, and $J_{D'}^\sim$ the N\'eron model over $D'$
of~$J_{D'}$.
\end{enumerate}
In particular, the morphism $f\colon X'\to J$ is proper if $X$ is
regular at all non-disconnecting double points of~$X_k$.
\end{proposition}
\begin{proof}
We apply the valuative criterion of properness. Arguments similar to
those in the proof of Proposition~\ref{prop9.1} show that conditions 1
and~3 are equivalent. So it remains to show that conditions 2 and~3
are equivalent.

We begin by giving a useful description of~$\Phi^\sim/\Phi$. Let $M$
denote the first homology group $\rH_1(G,\ZZ)$ of the graph~$G$; note
that it is canonically isomorphic to $\rH_1(G^\sim,\ZZ)$. We have:
\begin{subeqn}\label{eqn9.6.1}
M=\ker(d_*\colon\ZZ^{S^\sim_1}\to\ZZ^{S^\sim_0}), \quad\mbox{and}\quad 
M^\vee=\Hom_\ZZ(M,\ZZ)=\coker(d^*\colon\ZZ^{S^\sim_0}\to\ZZ^{S^\sim_1}). 
\end{subeqn}
Section~1 of \cite{Edixhoven2} shows that:
\begin{subeqn}\label{eqn9.6.2}
\Phi^\sim=\ZZ^{S^\sim_1}/(\ker(d_*)+\im(d^*))=
\ZZ^{S^\sim_0}[+]/\im(d_*d^*). 
\end{subeqn}
Raynaud's description of $\Phi$ in Section~\ref{section2} shows that
$\Phi$ is the image in $\Phi^\sim$ of the subgroup $\ZZ^{S_0}[+]$
of~$\ZZ^{S^\sim_0}[+]$. The inverse image under $d_*$ of this subgroup
is the relative homology group $\rH_1(G^\sim,S_0,\ZZ)$. Diagram~1.12
of \cite{Edixhoven2} says that the map from $\ZZ^{S^\sim_1}$ to
$\Phi^\sim$ factors through the map $\ZZ^{S^\sim_1}\to\ZZ^{S_1}$ that
sends $x$ in $S^\sim_1$ to $R(x)$ (recall that $R\colon X^\sim\to X$
is the resolution morphism). The image of $\rH_1(G^\sim,S_0,\ZZ)$ in
$\ZZ^{S_1}$ is the set of elements $a$ with $a(x)$ a multiple of
$e(x)$ for all $x$ in~$S_1$.  Let $q$ denote the quotient map from
$\ZZ^{S_1}$ to $\oplus_x\ZZ/e(x)\ZZ$.  Then we have:
\begin{subeqn}\label{eqn9.6.3}
\Phi^\sim/\Phi = (\oplus_x\ZZ/e(x)\ZZ)/qd^*\ZZ^{S_0}. 
\end{subeqn}
\begin{sublemma}\label{lemma9.6.4}
Let $Q$ be in $X(K)$, specializing to a non-disconnecting double point
$y$ of~$X_k$. Let $i$ be the integer such that $Q$ specializes to the
$i$th irreducible component in the chain of projective lines replacing
$y$ in~$X^\sim$, using the orientation of $G$ as in
\cite[2.2]{Edixhoven2}.  Then the image of $f(Q)$ in $\Phi^\sim/\Phi$
is represented by the element $\bar{i}$ in the factor at $y$ of
$\oplus_x\ZZ/e(x)\ZZ$.
\end{sublemma}
\begin{proof}
This follows directly from \cite[2.3]{Edixhoven2}. 
\end{proof}
Lemma~\ref{lemma9.6.4} implies that condition~3 above is equivalent to
the condition that, for every non-disconnecting double point $y$ of
$X_k$, the map from $\ZZ/e(y)\ZZ$ to $\Phi^\sim/\Phi$ is
injective. This injectivity is equivalent to the injectivity after
tensoring with $\ZZ_p$, for all~$p$, hence also to surjectivity in the
other direction after applying $\Hom({\cdot},\QQ_p/\ZZ_p)$, for
all~$p$. It follows that the map from $\ZZ/e(y)\ZZ$ to
$\Phi^\sim/\Phi$ is injective if and only if for every $p$ there
exists an element $g$ in $\rH_1(G,\QQ_p/\ZZ_p)$ such that $g(y)$ is
the class of $e(y)^{-1}$ and has $g(x)$ in $e(x)^{-1}\ZZ_p/\ZZ_p$ for
all $x$ in~$S_1$. Using that for any abelian group $A$ one has
$\rH_1(G,A)=A\otimes\rH_1(G,\ZZ)$, the existence of such a $g$ is seen
to be equivalent to the existence of a loop in $G$, containing $y$,
such that $v_p(e(y))\leq v_p(e(x))$ for all $x$ in that loop. This
finishes the proof of Proposition~\ref{prop9.6}.
\end{proof}
\begin{proposition}\label{prop9.7}
Suppose that $X_K$ has genus at least one.  The morphism $f\colon
X'\to J$ is a closed immersion if and only if the following two
conditions are satisfied:
\begin{enumerate}
\item every irreducible component of $X_k$ that is isomorphic to
$\PP^1_k$ contains a non-disconnec\-ting double point;
\item for every prime number $p$, every non-disconnec\-ting double
point $x$ of $X_k$ is contained in a loop in the graph $G$ such that
$v_p(e(x))\leq v_p(e(y))$ for every $y$ in that loop.
\end{enumerate}
\end{proposition}
\begin{proof}
Apply Propositions~\ref{prop9.5} and~\ref{prop9.6}. 
\end{proof}
\begin{corollary}\label{cor9.8}
Let $X$ be a regular proper nodal curve over a discrete valuation ring
$D$, with smooth geometrically irreducible generic fibre $X_K$ of
non-zero genus, with a given $P$ in~$X(K)$. Let $f_K\colon X_K\to J_K$
be the usual closed immersion sending $P$ to zero. Let $f\colon X^0\to
J$ be the induced morphism from the complement in $X$ of the set of
singular points of the special fibre $X_k$ of $X$ to the N\'eron model
over $D$ of~$J_K$. Then $f$ is a closed immersion if and only if all
double points of $X_{\bar{k}}$ are non-disconnecting, with $\bar{k}$
some algebraic closure of $k$.
\end{corollary}

\vspace{2\baselineskip}\noindent
{\bf Acknowledgements.} I would like to thank Torsten Ekedahl for 
pointing out to me the relation with log structures, as given in 
Section~\ref{section6}.

\vfill
\noindent
Bas Edixhoven\\
IRMAR\\
Campus de Beaulieu\\
35042 Rennes cedex\\
France

\end{document}